\documentclass[opre,nonblindrev,copyedit]{informs2_TR}

\usepackage{appendix}
\usepackage{endnotes}
\renewcommand{\v}[1]{\ensuremath{\boldsymbol{\mathrm{#1}}}}
\long\def\old#1{}
\usepackage{graphicx}
\usepackage{subfigure}

\usepackage{natbib}
 \bibpunct[, ]{(}{)}{,}{a}{}{,}%
 \def\bibfont{\small}%
\newcommand{\E}{\mathbb{E}}

\usepackage{multirow}
\usepackage{color}

\def\jnt#1{{\color{blue}#1}}

\long\def\old#1{}

\TheoremsNumberedThrough     
\ECRepeatTheorems

\EquationsNumberedThrough    

\begin{document}
\RUNAUTHOR{Tsitsiklis and Xu}

\RUNTITLE{Pricing of Fluctuations in Electricity Markets}

\TITLE{Pricing of Fluctuations in Electricity Markets}

\ARTICLEAUTHORS{%
\AUTHOR{John N. Tsitsiklis and Yunjian Xu} \AFF{Laboratory or
Information and Decision Systems, MIT, Cambridge, MA, 02139,
\EMAIL\{jnt@mit.edu, yunjian@mit.edu\}} }

\ABSTRACT{In an electric power system, demand fluctuations  may
result in significant ancillary cost to suppliers. Furthermore, in
the near future, deep penetration of volatile renewable electricity
generation is expected to exacerbate the variability of demand on
conventional thermal generating units. We address this issue by
explicitly modeling the ancillary cost associated with {demand
variability.}
 We {argue} that a time-varying price {equal to} the
suppliers' instantaneous marginal cost may not achieve social
optimality, and that consumer demand fluctuations should be
{properly} priced. We propose a dynamic pricing mechanism that
explicitly encourages consumers to adapt their consumption so as to
offset the {variability of} {demand on conventional units.}
 Through a dynamic game-theoretic
formulation, we show that (under suitable convexity assumptions) the proposed pricing mechanism achieves
 social optimality asymptotically, as the number of consumers increases
to infinity. Numerical results demonstrate that compared with
 marginal cost pricing, the proposed mechanism creates a stronger incentive for consumers to shift their peak
load, and therefore has the potential to reduce the need for
long-term investment in peaking plants.}
%


\KEYWORDS{Electricity market, dynamic pricing, social welfare\\
{{\it Date:} \today}}

%
%


\maketitle



%

\section{Introduction}\label{sec:intro}

{This paper is motivated by} the fact that fluctuations in the
demand {for electricity} to be met by conventional thermal
generating units typically result in significantly increased, and
nontrivial,  ancillary costs. Today, such demand fluctuations are
mainly due to time-dependent
  consumer preferences. In addition, in the future, a certain percentage of electricity
production is required by law in many states in the U.S. to come
from renewable resources \citep{BW08}. The dramatic
 volatility  of renewable energy
 resources may aggravate the variability of the demand for conventional thermal
generators and result in significant ancillary cost.
 More concretely, either a demand surge or a decrease in renewable generation may result in
(i) higher energy costs due to the deployment of peaking plants with higher ramping rates but higher marginal cost, such as
 oil/gas combustion turbines, and (ii) the redispatch cost\footnote{A
certain level of reserve must always be maintained in an electric
power system. Local reserve shortages are usually due to the quick
increase of system load rather than a capacity deficiency. If the
increase of system load makes the system short in reserves, the
system will redispatch resources to increase the amount of reserves
available. Redispatch generally increases the generation cost
and results in higher prices. The redispatch cost
can be very high (cf. Section 2.3.2 of \cite{E10}).} that the system will incur to meet reserve
constraints if the increase of demand (or decrease of renewable generation) causes a reserve shortage.

There is general agreement that
  charging real-time prices (that reflect current operating conditions)
 to electricity consumers has the potential of reducing supplier ancillary cost, improving
 system efficiency, and lowering volatility in wholesale prices \citep{DOE06,SL08,C10}.
Therefore, dynamic pricing, especially real-time marginal cost
pricing, is often identified as a priority for the implementation of
 wholesale electricity markets with responsive demand \citep{H10},
 which in turn raises many new questions. For example,
 should prices for a given time interval be calculated ex ante or ex post? Does real-time pricing introduce
 the potential for new types of market instabilities? How is supplier competition affected? In this paper, we abstract away
 from {almost all of these} questions and focus on the specific issue of whether prices should also explicitly encourage consumers to adapt their
demands so as to reduce supplier ancillary cost.

\begin{figure}
\begin{minipage}[t]{0.5\linewidth}
\includegraphics[width=3.02in]{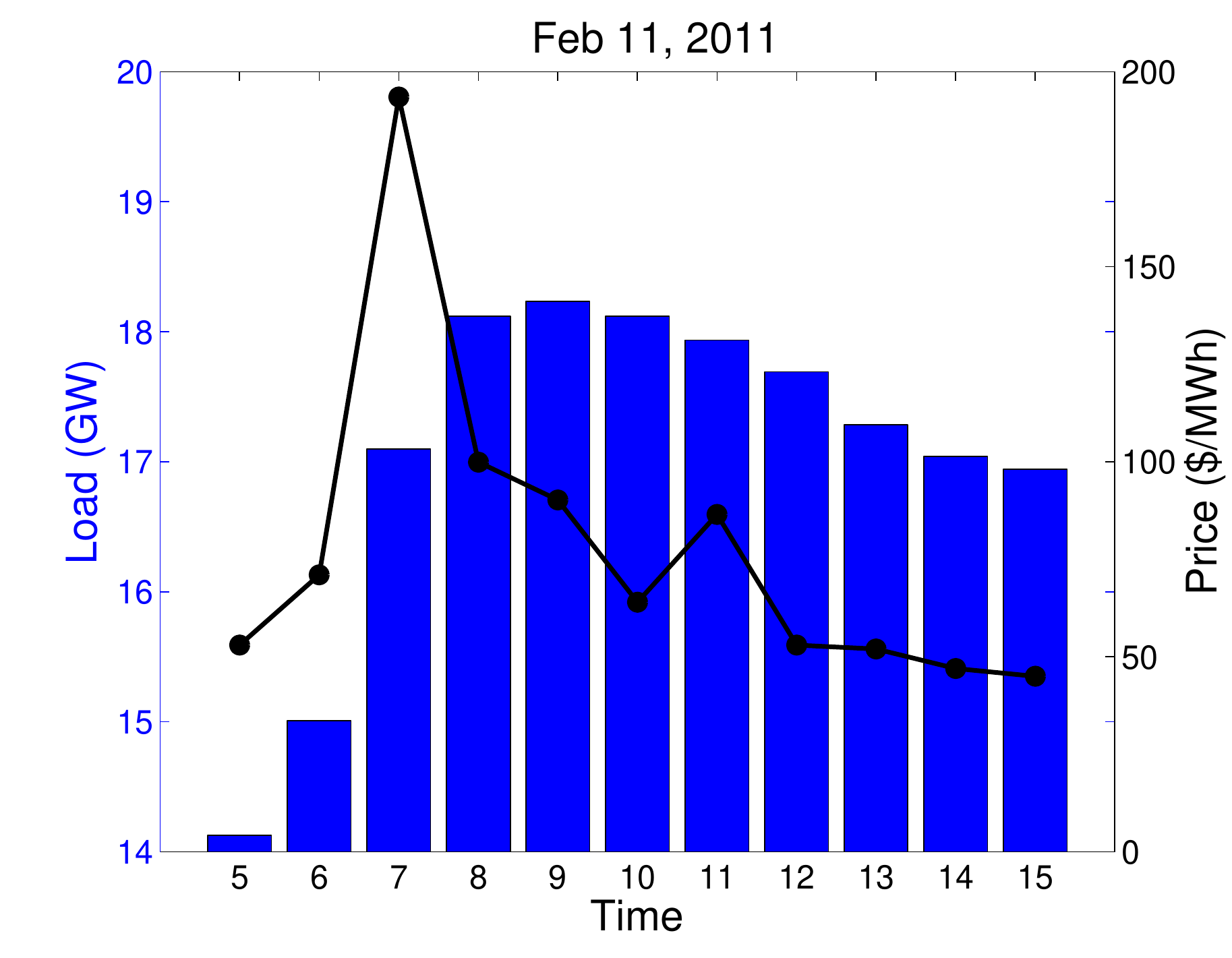}
\end{minipage}%
\begin{minipage}[t]{0.5\linewidth}
\centering
\includegraphics[width=3.01in]{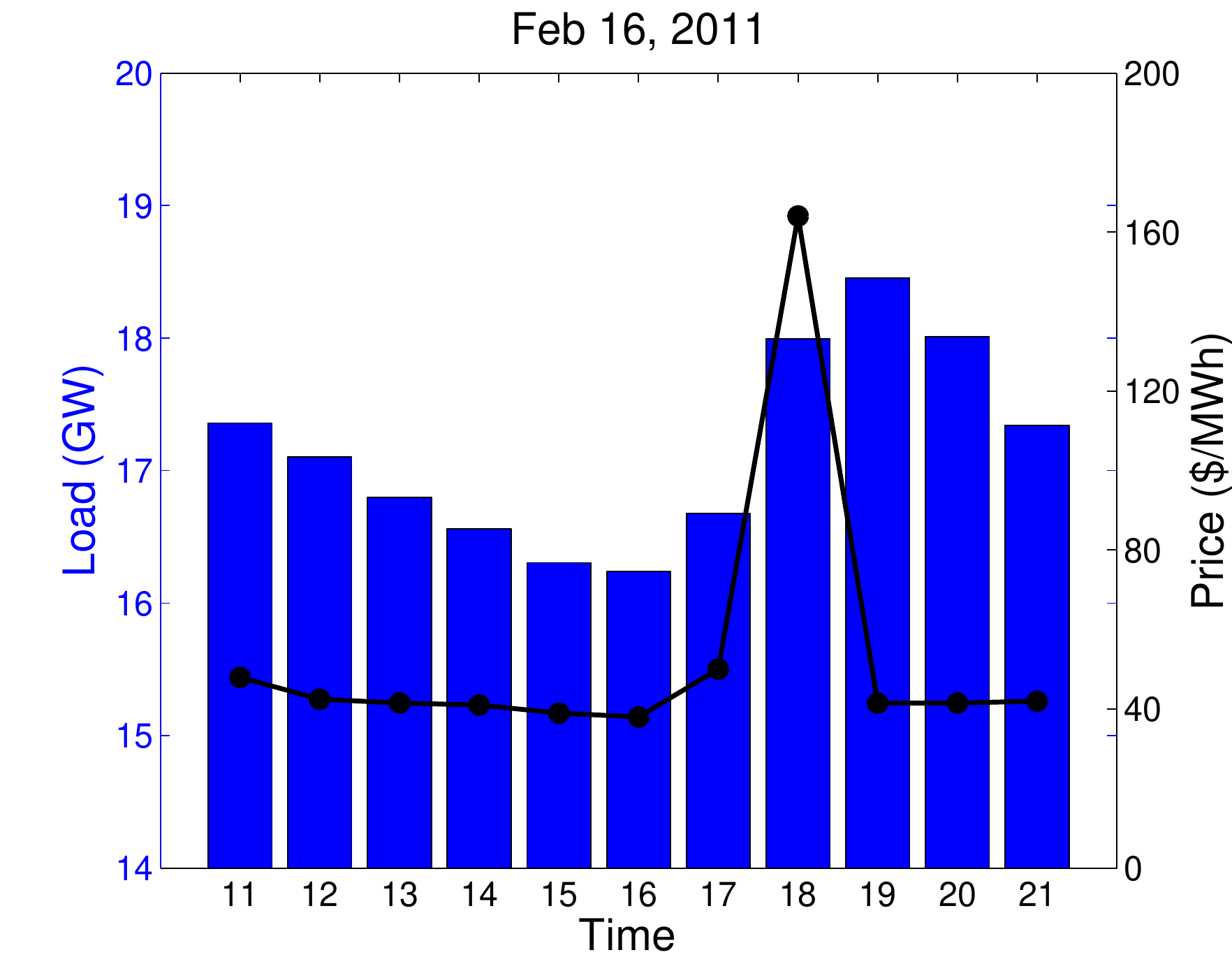}
\end{minipage}
\caption{Real-time prices and actual system load, ISO New England
Inc. Blue bars represent the real-time system loads and the dots
connected by a black line represent the hourly
prices.}\label{Fig:211}
\end{figure}

To illustrate the issue that we focus on, we note that a basic model
of electricity markets assumes that the cost of satisfying a given
level $A_t$ of aggregate demand during period $t$ is of the form
$C(A_t)$. It then follows that in a well-functioning wholesale
market, the observed price should more or less reflect the marginal
cost $C'(A_t)$. {In particular,} prices should be more or less
determined by the aggregate demand level. Empirical data do not
quite support this view. Fig. \ref{Fig:211} plots the real-time
system load and the hourly prices on February 11, 2011 and on
February 16, 2011, as reported by the New England ISO \citep{F11}.
We observe that prices do not seem to be determined solely by $A_t$
{but that the} changes in demand, $A_t-A_{t-1}$, {also} play a major
role. In particular, the largest prices seem to occur after a demand
surge, and not necessarily at the hour when the load is highest.
We take this as evidence that the total cost over $T+1$ periods is
not of the form
$$
\sum_{t=0}^T C(A_t),
$$
but rather of the form
\begin{equation}\label{eq:aa}
\sum_{t=0}^T \big(C(A_t)+H(A_{t-1},A_{t})\big),
\end{equation}
for a suitable function $H$.

We take the form of Eq.\ (\ref{eq:aa}) as our starting point and raise the question of the appropriate prices.
A naive view would argue that at time $t$, $A_{t-1}$ has already been realized, and taking its value for granted,
a consumer should be charged a unit price equal to
\begin{equation}\label{eq:two}
C'(A_t) +\frac{\partial}{\partial A_t} H(A_{t-1},A_t),
\end{equation}
which is the supplier's marginal cost at stage $t$. We refer to this
naive approach as ``marginal cost pricing'' (MCP). However, a simple
argument based on standard mathematical programming optimality
conditions shows that for system optimality to obtain, the demand
$A_{t-1}$ should also incur
 (after $A_t$ is realized) a unit price of
\begin{equation}\label{eq:three}
\frac{\partial}{\partial A_{t-1}} H(A_{t-1},A_t).
\end{equation}
In day-ahead markets, suppliers typically carry out an intertemporal
optimization, and it is reasonable to expect that the two types of
marginal costs, captured by Eqs.\ \eqref{eq:two} and
\eqref{eq:three}, are both properly accounted for. However, in
current real-time balancing markets, once $A_{t-1}$ is realized, a
supplier will aim at charging the marginal cost in Eq.\
\eqref{eq:two}, but will be unable to charge the additional marginal
cost in Eq.\ \eqref{eq:three} to the past demand $A_{t-1}$. In
contrast, the pricing mechanism that we propose and analyze in this
paper is designed to include the additional marginal cost in Eq.
\eqref{eq:three}.\footnote{In current two-settlement systems, the
real-time prices are charged only on the difference of the actual
demand and the estimated demand at the day-ahead market. However,
the two-settlement  system provides the same real-time incentives to
price-taking consumers, as if they were purchasing all of their
electricity at the real-time prices (cf. Chapter 3-2 of
\cite{S02}).}

The actual model that we consider will be richer from the one discussed above
 in a number of respects. It includes an exogenous source of uncertainty (e.g., {representing} 
 weather conditions) that has an impact on consumer utility
  and supplier cost, {and therefore the model can incorporate the effects of volatile} renewable electricity production.
    It  allows for consumers with internal state variables (e.g., a consumer's
  demand may be affected by how much electricity she has already used). It also allows for multiple consumer types
  (i.e., with different utility functions and different internal state dynamics). Consumers are generally modeled as price-takers,
  as would be the case in a model involving an infinity (a continuum) of consumers. However, we also consider the case of finite
  consumer populations and explore certain equilibrium concepts that are well-suited to the case of finite but large consumer populations.
  On the other hand, we ignore most of the distinctions between ex post and ex ante prices. Instead, we assume that at each time step,
  the electricity market clears. The details of how this could happen are important, but are generic to electricity markets,
  hence not specific to our models, and somewhat
  orthogonal to the subject of this paper. (See however Appendix \ref{sec:im} for some discussion of implementation issues.)

The ancillary cost function $H(A_{t-1},A_t)$ is of course a central
element of our model. How can we be sure that this is the right
form? In general, redispatch and reserve dynamics are complicated
{and} one should not expect such a function to capture all of the
complexity of the true system costs; perhaps, a more complex
functional form such as $H(A_{t-2},A_{t-1},A_t)$ would be more
appropriate. We believe that the form we have chosen is a good
enough  approximation, at least under certain conditions. To argue
this point, we present in Appendix \ref{sec:A0} an example that
involves a more detailed system model (in which the true cost is a
complicated function of the entire history of demands) and show that
a function of the form $H(A_{t-1},A_t)$
 can capture most of the cost of ancillary services.

\subsection{Summary of contributions}
Before continuing, we provide here a roadmap of the paper together
with a summary of our main contributions.
\begin{itemize}{\itemsep=0pt}
\item[(a)]
We provide a stylized (yet quite rich) model of an electricity
market, which  incorporates the cost of ancillary services (cf. Section \ref{sec:model}).
\item[(b)] We provide some justification of the form of the cost function in our model, as a reasonable approximation of more detailed physical models
(cf. Appendix \ref{sec:A0}).
\item[(c)] We propose a pricing mechanism that properly charges for the effects of consumer actions on  ancillary services (cf. Section \ref{sec:PM}).
\item [(d)] For a continuum model involving nonatomic price-taking consumers, we consider  Dynamic Oblivious Equilibria (DOE),
in which every consumer maximizes her expected payoff under the
sequence of prices induced by {a} DOE strategy profile (Section
\ref{sec:DOE}). We show that (under standard convexity assumptions),
our mechanism maximizes social welfare (cf. Theorem \ref{thm:SO_social} {in Section \ref{sec:social}}).

\item[(e)] We {carry out a game-theoretic analysis of} the case of a large but finite number of consumers.
 We show that a large population of consumers who act
according to a DOE (derived from an associated continuum game)
 results in asymptotically optimal (as the number of consumers goes to infinity) social welfare (cf. Theorem \ref{thm:SO_social}
in Section \ref{sec:social}), {and asymptotically maximizes every
consumer's expected payoff} (this is an ``asymptotic Markov
equilibrium'' property; cf.\ Theorem \ref{thm:APE} {in Section
\ref{sec:AMP}}).

\item[(f)] We illustrate the potential benefits of our mechanism through a simple numerical example. In particular, we show that compared with
marginal cost pricing, the proposed mechanism reduces the peak load,
and therefore has the potential to reduce the need for long-term
investments in peaking plants (cf.\ Appendix \ref{sec:num}).

\end{itemize}

\subsection{Related literature}

There are two streams of literature, {on electricity pricing and on
game theory, that are relevant to our work, and which we now proceed
to discuss, while also highlighting the differences from the present
work.}

Regarding electricity markets, the impact of supply friction on
economic efficiency and price volatility has received some recent
attention.  \cite{M08} shows that under ramping constraints, the
prices faced by consumers may not necessarily equal the true
supplier marginal cost. In a continuous-time competitive market
model, \cite{CM10} show that the limited capability of generating
units to meet real-time demand, due to
 relatively low ramping rates, does not harm social welfare,
but may result in extreme price fluctuations. In a similar spirit,
\cite{KM10} construct a dynamic game-theoretic model to study the
tradeoff between economic efficiency and price volatility. Closer to
the present paper, \cite{CM10b} construct a dynamic newsboy model to
study the reserve management problem in electricity markets, where
the demand is assumed to be exogenous. The supplier cost in their
model depends not only on the overall demand, but also on the
generation resources used to satisfy the demand. For example, a
quickly increasing demand may require more responsive and more
expensive resources (e.g., peaking generation plants).

{To study the impact of pricing mechanisms on consumer behavior and
load fluctuations, we construct a dynamic game-theoretic model that
differs from existing dynamic models for electricity markets and
incorporates both the consumers' responses to real-time price
fluctuations and the suppliers' ancillary cost incurred by load
swings.}
Some major differences between our model and existing ones are discussed
at the end of Section \ref{sec:model}.

On the game-theoretic side,  {the standard solution concept for
stochastic dynamic  games is {the} \textbf{Markov perfect
equilibrium} (MPE) \citep{FT91,MT88}, {involving strategies} where
an agent's action depends on the current state of all agents. As the
number of agents grows large, the computation of an MPE is often
intractable \citep{DP07}. For this reason, alternative equilibrium
concepts, for related games featuring a nonatomic continuum of
agents (e.g., ``oblivious equilibrium'' or ``stationary
equilibrium'' for dynamic games without aggregate shocks), have
received much recent attention \citep{WBV08,AJ10}.}

There is a large literature on a variety of approximation properties
of nonatomic equilibria \citep{MV93,A04,A08}. Recently,
\cite{AJ10} derive
 sufficient conditions for a
stationary equilibrium strategy to have the \textbf{Asymptotic
Markov Equilibrium} (AME) property, i.e., a stationary
equilibrium strategy asymptotically maximizes every agent's expected
payoff (given that all the other agents use the same stationary
equilibrium  strategy), as the number of agents grows large. {Their
model includes random shocks that} are assumed to be idiosyncratic
across agents. However, in the problem that we are interested in, it
is important to incorporate aggregate shocks (such as weather
conditions) that have a global impact on all agents. In this spirit,
\cite{WBV10} consider a market model with aggregate profit shocks,
and study an equilibrium concept at which every firm's strategy
depends on the firm's current state and on the recent history of the
aggregate shock. For a general dynamic game model with aggregate
shocks, \cite{BC10} shows that a nonatomic {counterpart} of an MPE,
which we refer to as a \textbf{Dynamic Oblivious Equilibrium (DOE)}
in this paper, asymptotically approximates an MPE in the sense that
as the number of agents increases to infinity, the actions taken in
an MPE can be well approximated by those taken by a DOE strategy of
the nonatomic limit game. However, without further restrictive
assumptions on the agents' state transition kernel, the
approximation property of the actions taken by a DOE strategy
 does not imply the AME property of the DOE, and we are not aware of
 any AME results for models that include aggregate shocks. Our work is different in this respect:
for a dynamic nonatomic model with aggregate shocks, which
is a simplified variation of the general model considered in
\cite{BC10}, we prove the AME property of a DOE.

The efficiency of nonatomic equilibria for static games has been
addressed in recent research \citep{RT04,M04,BC11}. For a dynamic
industry model with a continuum of identical producers and exogenous
aggregate shocks, \cite{LP71} show (under convexity assumptions)
that the expected social welfare is maximized at a unique
competitive equilibrium. In a similar spirit, in this paper we show
(under convexity assumptions) that the proposed pricing mechanism
maximizes the expected social welfare in a model involving a
continuum of (possibly heterogeneous) consumers. We also consider
the case of a large but finite number of consumers, and show that
the expected social welfare is approximately maximized if all
consumers act according to a nonatomic equilibrium (DOE). For large
dynamic games, the asymptotic social optimality of nonatomic
equilibria (DOEs) established
 in this paper seems to be new.

\old{
\subsection{Outline of the paper}

The rest of the paper is organized as follows. In Section
\ref{sec:model}, we  introduce the dynamic game model. In Section
\ref{sec:PM}, inspired by the observation that the marginal cost
pricing may not be socially optimal in a market with friction (see
Example \ref{ex:price}), we propose an alternative dynamic pricing
mechanism. In Section \ref{sec:DOE}, we give the formal definition
of a DOE. In Section \ref{sec:AMP}, we show that as the number of
consumers increases to infinity, every consumer's expected payoff is
approximately maximized by a DOE strategy, if the other consumers
follow the DOE strategy. In Section \ref{sec:social}, we show that
under the proposed pricing mechanism, social welfare is maximized at
a DOE of a continuum game. Also, for a dynamic game with finitely
many consumers, and under a certain convexity assumption, we
validate the asymptotic social optimality of the proposed pricing
mechanism: as the number of consumers grows large, if all consumers
use a DOE strategy, then social welfare is approximately maximized.
In Section \ref{sec:num}, we present several numerical examples to
compare the proposed pricing mechanism with marginal cost pricing.
Finally, in Section \ref{sec:con}, we make some brief concluding
remarks and discuss some directions for future work.}

\section{Model}\label{sec:model}

We consider a $(T+1)$-stage dynamic game with the following
elements:

\begin{enumerate}{\itemsep=0pt}

\item The game is played in \textbf{discrete time}. We index the time
periods with
 $t=0,1,\ldots,T$. Each stage may represent a five minute interval in real-time balancing markets
 where prices and dispatch solutions are typically provided at five minute intervals.

\item There are $n$ \textbf{consumers}, indexed by $1,\ldots,n$.

\item At each stage $t$, let $s_t \in \mathcal {S} $ be an \textbf{exogenous
state}, which evolves as a Markov chain and whose transitions are
{not affected by} consumer actions. The set $\mathcal {S} $ is
assumed to be finite. In electricity markets, the exogenous  state
may represent time and/or weather conditions, which impact
consumer utility and supplier cost. It may also represent the level
of renewable generation.

\item  For notational conciseness, for $t \ge 1$, let
$\overline s_t=(s_{t-1},s_t)$, and let $\overline s_0=s_0$. We use
$\overline {\mathcal {S}}_t$ to denote the set of all possible
$\overline s_t$.  We  refer to $\overline s_t$ as the \textbf{global
state} at stage $t$.

\item Given an initial global state $s_0$,
the \textbf{initial states (types) of the consumers},
$\{x_{i,0}\}_{i=1}^n$,
 are independently drawn according to a
probability measure $\eta_{s_0}$ over a finite set $\mathcal {X}_0$.
We use $X$ to denote the cardinality of $\mathcal {X}_0$.

\item
The \textbf{state of consumer} $i$ at stage $t$ is denoted by
$x_{i,t}$. At $t=0$, consumer $i$'s initial state, $x_{i,0}$,
indicates her type. For $t=1,\ldots,T$, we have
$x_{i,t}=(x_{i,0},z_{i,t})$, where $z_{i,t} \in \mathcal {Z}$ and
$\mathcal {Z}=[0,Z]$ is a compact subset of $\mathbb{R}$. The
variables $\{z_{i,t}\}_{i=1}^n$ allow us to model  intertemporal
substitution effects
 in consumer $i$'s demand.

\item We use $\mathcal {X}_t$ to denote a \textbf{consumer's state space} at
stage $t$. In particular, at stage $t \ge 1$, $\mathcal
{X}_t=\mathcal {X}_0 \times \mathcal {Z}$.

\item At stage $t$,  consumer $i$ takes an \textbf{action} $a_{i,t}$
and receives a nonnegative \textbf{utility}\footnote{At $t = 0$,
$U_0$ is a mapping  from $\mathcal {X}_0  \times \mathcal {S} \times
\mathcal {A}$ to $[0,\infty) $, while for $t \ge 1$, $U_t$ is a
mapping from $\mathcal {X}_0 \times {\mathcal {Z}} \times \mathcal
{S} \times \mathcal {A}$ to $[0,\infty) $.}
$U_t(x_{i,t},s_t,a_{i,t})$.

\item Each consumer's \textbf{action space} is
$\mathcal {A}=[0,B]$, where $B$ is a positive real number. (In the electric power context, $B$ could reflect
 a local transmission capacity constraint.)

\item  We
 use $A_t=\sum\nolimits_{i=1}^n a_{i,t}$ to denote the \textbf{aggregate
 demand} at stage $t$.

\item Given consumer $i$'s
 current state, $x_{i,t}$, and the next exogenous state $s_{t+1}$,
the
 next state of consumer $i$ is determined by
her action
 taken at stage $t$, i.e., $x_{i,t+1}=\left( x_{i,0},
z_{i,t+1}\right)$, where $z_{i,t+1}=
 r(x_{i,t},a_{i,t},s_{t+1})$, for a given function $r$.

\item Let $G_t=A_t+R_t$ be the \textbf{capacity} available at stage $t$,
where $R_t$ is the system \textbf{reserve} at stage $t$. For
simplicity, we assume that the system reserve at stage $t$ depends
only on the current aggregate demand, $A_t$, and the current
exogenous state $s_t$. That is, we have $R_t= g(A_t,s_t)$ for a
given function of $g$ that reflects the reserve policy of the system
operator.

\item At stage $t$, let $\overline C(A_t, R_t, s_t)$ {be the total conventional \textbf{generation cost}, that is, the sum}
of the supplier's cost to meet the aggregate
demand $A_t$ through its primary energy resources, e.g., base-load
power plants, and the cost to maintain a system reserve $R_t$.
 Since $R_t$ depends only on $A_t$ and $s_t$, we can write
$\overline C(A_t, R_t, s_t)$ as a function of $A_t$ and $s_t$, i.e.,
there exists a \textbf{primary cost function} $C: \mathbb{R} \times
\mathcal {S} \to [0,\infty)$ such that $C(A_t,s_t)=\overline C(A_t,
R_t, s_t)$. We assume that for any $s \in \mathcal {S}$,
$C(\cdot,s)$ is nondecreasing.

\item  At stage $t \ge 1$, let $\overline H(A_{t-1},A_{t},R_{t-1},R_{t},s_t)$ denote
the ancillary
 cost incurred by load swings\footnote{In general, the supplier
ancillary cost may depend on the entire history of system load and
global states. However, {ancillary cost functions with the simple
form} $\overline H(A_{t-1},A_{t},R_{t-1},R_{t},s_t)$ can serve as a
good approximation of the supplier's true ancillary cost (cf.
Appendix \ref{sec:A0}). }. Since $R_t$ depends only on $A_t$ and
$s_t$, we can write $\overline H(A_{t-1},A_{t},R_{t-1},R_{t},s_t)$
as a function of $A_{t-1}$, $A_t$, $s_{t-1}$, and $s_t$, i.e., there
exists an \textbf{ancillary cost function} $H: \mathbb{R}^2 \times
\mathcal {S}^2 \to [0,\infty)$ such that $\overline H = H(A_{t-1},
A_t, \overline s_t)$. The ancillary cost at stage $0$ is assumed to
be a function of $s_0$ and $A_0$.

\item At stage $0$,  the total supplier cost is of the form
\begin{equation}\label{equa:cost0}
C(A_0,s_0)+  H_0(A_{0}, s_0),
\end{equation}
and for  $t=1,\ldots,T$,
 the total supplier cost at
stage $t$ is given by
\begin{equation}\label{equa:cost}
C(A_t,s_t)+   H(A_{t-1}, A_t, \overline s_t).
\end{equation}

\end{enumerate}


In contrast to existing dynamic models for electricity markets with
an exogenous demand process \citep{CM10,CM10b}, our dynamic
game-theoretic model incorporates the consumer reactions to price
fluctuations, and allows us to study the impact of pricing
mechanisms on consumer behavior and economic efficiency. Through a
dynamic game-theoretic formulation, \cite{KM10} study the tradeoff
between economic efficiency and price volatility. Our model is
different from the one studied in \cite{KM10} in the following
{respects:}

\begin{enumerate}{\itemsep=0pt}
\item
 Our model allows the generation cost to
depend on  an exogenous state, and therefore {can incorporate}
supply-side volatility due to  uncertainty in renewable electricity
generation. As an example, consider a case where the exogenous
state, $s_t$, represents the electricity generation from renewable
resources at stage $t$. Then the demand for conventional generation
is $A_t-s_t$. Suppose that the system reserve is proportional to the
system load, say, $\delta A_t$ {for some constant} $\delta>0$. The
cost function, $\overline C(A_t,R_t,s_t)$, then depends only on the
output of conventional generating units, $A_t-s_t$, and the system
reserve, $\delta A_t$. The ancillary cost occurred at stage $t$
depends on the system reserve and the outputs of conventional
generating units at stages $t-1$ and $t$, and is therefore a
function of $A_{t-1}$, $A_t$, $s_{t-1}$, and $s_t$.

\item More important, instead of penalizing
each consumer's attempt to change her {own} action across time,
 the ancillary cost function in our model penalizes the change in
 the aggregate demand by all consumers.  The change in a single
 consumer's action may harm or benefit the social welfare, while the
 volatility of the aggregate demand is usually undesirable.
\end{enumerate}

The main feature of our model is the ancillary cost function $H$,
which makes the supplier cost nonseparable over time. In an electric
power system, the ancillary cost function models the costs
associated with the variability of conventional thermal generator
output, such as the energy cost of peaking plants and the redispatch
cost. Note that the ancillary cost is not necessarily zero when $A_t
\le A_{t-1}$, because thermal generating units have ramping-down
constraints, and because a decrease in renewable electricity
production may {lead to an increase of} the system reserve, even if
$A_t \le A_{t-1}$.
 The presence of the ancillary cost function makes {conventional}
 marginal cost pricing  inefficient (cf. Example \ref{ex:price} in Section \ref{sec:PM}).

\old{ For the energy conservation problem in data centers, A similar
model is used to study the energy consumption in data centers in
\cite{LW11}, where the switching cost}

To keep the model simple, we do not incorporate any idiosyncratic
randomness in  consumer state evolution. Thus, besides the
randomness of consumer types (initial states),
 the only source of stochasticity in the model {is}
the exogenous state $s_t$.

To effectively highlight the impact of pricing mechanisms on
consumer behavior, as well as on economic efficiency and demand
volatility, we have made the following simplifications and
assumptions for the power grid:

\begin{itemize}
\item [(a)]
 As in \cite{CM10}, we assume that the physical production capacity is  large enough so that the possible
 changes
of the generation capacity are not constrained.

\item [(b)]  Transmission capacity is large enough to avoid any
congestion. We also assume that the cost of supplying electricity to
consumers {at} different locations is the same. Therefore, a common
price for all consumers is appropriate.

\item [(c)]  We use a simplified form of ancillary cost functions,
$\overline H(A_{t-1},A_{t},R_{t-1},R_{t},s_t)$, to approximate the
supplier ancillary cost. \old{ Such an
approximation is well justified if the ancillary energy service,
e.g., a peaking plant, is deployed only at a single stage. At
the next stage, either the demand drops, or the primary
energy resource can ramp up rapidly enough to meet the increasing
demand. This is usually the case, even in current electricity
markets where most consumers pay a fixed price for the electricity
they consume. For example, among the $239$ deployments of New York
independent system operator reserve in $2002$, most of them were
$12$ minutes or less in duration \cite{K03}. Under real-time
pricing, it may not be implausible to suppose that this
approximation will be even better justified with the price-responsive
demand.}In Appendix \ref{sec:A0}, we present a numerical example to
justify this approximation.

\old{ Numerical results show that the difference approximation

In general, the supplier ancillary cost may depend on the entire
history of system load and exogenous states. We assume that the system
reserve is a function of the current aggregate demand and the
current exogenous state. Under this assumption,

This assumption allows us to consider a simplified form of ancillary
cost functions, which depends only on $A_{t-1}$, $A_t$ and
$\overline s_t$.}
\end{itemize}

\section{The Pricing Mechanism}\label{sec:PM}

The marginal cost pricing mechanism discussed in Section \ref{sec:intro}
charges a time-varying unit
price on each consumer's demand.
As demonstrated in the following example, a time-varying price that equals the supplier's instantaneous
  marginal cost may not achieve social optimality in a {setting that includes ancillary costs.} For this reason, we
propose a new pricing mechanism that takes into account the
ancillary cost associated with a consumer's demand at the {{\bf
previous}} stage.

\begin{example}\label{ex:price}  \rm{
Consider a two-stage deterministic model with one
 consumer and one supplier. At stage $t$, the consumer's utility function
is $U_t:[0,\infty) \to [0,\infty)$.
  Let
$a_t$ denote the demand at stage $t$, and let $\v{a}=(a_0,a_1)$.
  Let
$g_t$ denote the actual generation at stage $t$, and let
$\v{g}=(g_0,g_1)$.  Two unit prices, $p_0$ and $p_1$, are charged on
the consumption at stage $0$ and $1$, respectively. Let
$\v{p}=(p_0,p_1)$.  The consumer's payoff-maximization problem is
\begin{equation}\label{eq:exa_C}
 { \mathop {\rm Maximize}\limits_{\v{a}}}\;\;\;\;\;
U_0(a_0)- p_0 a_0 + U_1(a_1)- p_1 a_1.
\end{equation}
Let $H_0$ be identically zero, and let the ancillary cost function
at stage $1$ depend only on the difference between the supply at the
two stages. That is, the ancillary cost at stage $1$ is of the form
$H(g_1-g_0)$. The supplier's profit-maximization problem is
\begin{equation}\label{eq:exa_S}
\displaystyle { \mathop {\rm Maximize}\limits_{\v{g}}}\;\;\;\;\; p_0
g_0 + p_1 g_1-C(g_0)-C(g_1)-H(g_1-g_0).
\end{equation}
The social planner's problem is
\begin{equation}\label{eq:exa_W}
\begin{array}{l}
\displaystyle {\mathop {\rm Maximize}\limits_{(\v{a},\v{g})}}\;\;\;\;\; U_0(a_0)+
U_1(a_1)-C(g_0)-C(g_1)-H(g_1-g_0) \\[5pt]
\displaystyle {\rm{subject}}\;{\rm{to}}\;\;\;\;\;\v{a} = \v{g}.
\end{array}
\end{equation}

Now consider a \textbf{competitive equilibrium},
$(\v{a},\v{g},\v{p})$, at which the vector $\v{a}$ solves the consumer's optimization
problem \eqref{eq:exa_C}, the vector $\v{g}$ solves the supplier's optimization
problem \eqref{eq:exa_S}, and the market clears, i.e., $\v{a}=\v{g}$.
Suppose that the utility functions are concave and continuously
differentiable, and that the cost functions $C$ and $H$ are convex and
continuously differentiable. We further assume that $H'(0)=0$,
 and that for $t=0,1$, $U'_t(0)>C'(0)$, $U'_t(B)<C'(B)$. Then, there exists a competitive
equilibrium, $(\v{a},\v{g},\v{p})$, which satisfies the following
conditions:
\begin{equation}\label{eq:exa_C1}
\left\{ \begin{array}{l}
\displaystyle U'_0({a_0}) = {p_0},\\[5pt]
\displaystyle U'_1({a_1}) = {p_1},
\end{array} \right. \;\;\;\;\;\;\;\;\;\;\;\;\;\;\;\;\;\;\left\{ \begin{array}{l}
\displaystyle C'({a_0}) -H'(a_1-a_0) = {p_0},\\[5pt]
\displaystyle C'({a_1}) +H'(a_1-a_0)= {p_1}.
\end{array} \right.
\end{equation}
We conclude that the competitive equilibrium solves the social
welfare maximization problem in \eqref{eq:exa_W}, because it satisfies the following
(sufficient) optimality conditions:
\begin{equation}\label{eq:exa_C2}
\begin{array}{l}
 \displaystyle U'_0({a_0})   = C'(a_0)-H'(a_1-a_0),\qquad  U'_1({a_1})  = C'(a_1)+ H'(a_1-a_0),\\[5pt]
 \displaystyle a_0=g_0,\;\;\;a_1=g_1.
 \end{array}
\end{equation}

However, we observe that the socially optimal price $p_0$ does not
equal the supplier's {instantaneous} marginal cost at stage $0$,
$C'(a_0)$. Hence, by setting the price equal to $C'(a_0)$, {as would
be done in a real-time balancing market,} we may not achieve social
optimality. More generally, marginal cost pricing need not be
socially optimal because it does not take into account the
externality conferred by the action $a_0$ on the ancillary cost at
stage $1$, $H(a_1-a_0)$. At a socially optimal competitive
equilibrium, the consumer should pay
$$
\big(C'(a_0)-H'(a_1-a_0)\big)a_0+\big(C'(a_1)+H'(a_1-a_0)\big)a_1,
$$
{i.e.,} the price on $a_0$ {should be} the sum of the supplier
marginal cost at stage $0$, $C'(a_0)$, and the marginal ancillary
cost associated with $a_0$, $-H'(a_1-a_0)$, {which is determined at
the next stage, after $a_1$ is realized}. } $\;\;\;\;\;\;\square$
\end{example}

Before {describing the precise} pricing mechanism {we propose}, we
introduce a differentiability assumption on the cost functions.

\begin{assumption}\label{A:price}
For any $s \in \mathcal {S}$, $C(\cdot,s)$ and $H_0(\cdot,s)$ are continuously
differentiable on ${[0},\infty)$.
 For any $(A',\overline s) \in \mathcal {A} \times \mathcal {S}^2$, $H(A,A',\overline s)$ and $H(A',A,\overline s)$
 are continuously differentiable in $A$ on $[0,\infty)$.\footnote{At the boundary of
the domain, $0$, we require continuity of the
right-derivatives of $C$, $H_0$, and $H$.}
 \end{assumption}

Inspired by Example \ref{ex:price}, we introduce prices
\begin{equation}\label{equa:pricep}
{p_t} =C'(A_t, s_t ),\qquad t=0,\ldots,T,
\end{equation}
and
\begin{equation}\label{equa:priceq}
q_t =  \dfrac{\partial H(A_{t-1}, A_t, \overline s_t )}{\partial
A_{t-1}},\quad w_t =  \dfrac{\partial H(A_{t-1}, A_t, \overline s_t
)}{\partial A_{t}},  \qquad t=1,\ldots,T.
\end{equation}
At stage $0$, we let $q_0=0$ and $w_0=H_0'(A_0,s_0)$.
 Under the proposed pricing mechanism, consumer $i$'s payoff at stage $t$ is
given by
\begin{equation}\label{equa:Spayoff1}
U_t(x_{i,t},s_t,a_{i,t})  - (p_t+w_t) a_{i,t} -   q_t a_{i,t-1} .
\end{equation}
Note that $p_t+w_t$ is the supplier marginal cost at stage $t$
(including the marginal ancillary cost). The proposed pricing
mechanism charges consumer $i$ an additional price $q_t$ {on her
previous demand}, equal to the marginal ancillary cost with respect
to $a_{i,t-1}$.

We now {define some of the notation that we will be using.}
 For $t=1,\ldots,T$, let
$y_{i,t}=(a_{i,t-1},x_{i,t})$ be the \textbf{augmented state} of
consumer $i$ at stage $t$. At $t=0$,
 let $y_{i,0}=x_{i,0}$.
{For} stage $t$, let $\mathcal {Y}_t$ be the set of all possible
augmented states. In particular, we have $\mathcal {Y}_0=\mathcal
{X}_0$, and $\mathcal {Y}_t=\mathcal {A} \times \mathcal {X}_t$, for
$t=1,\ldots,T$.

{Let $\Delta_n(D)$ be} the set of empirical probability
distributions over a given set {$D$} that can be generated by $n$
samples from $D$. (Note that empirical distributions are always
discrete, even if $D$ is a continuous set.) Let $f_{t} \in
\Delta_n(\mathcal {Y}_t)$ {be} the empirical distribution of the
augmented state of all consumers at stage $t$, and {let} $f_{-i,t}
\in \Delta_{n-1}(\mathcal {Y}_t)$ {be} the empirical distribution of
the augmented state of all consumers (excluding consumer $i$) at
stage $t$. {We} refer to $f_t$ as the \textbf{population state} at
stage $t$. Let $u_{t} \in \Delta_n(\mathcal{A})$ denote the
empirical distribution of all consumers' actions at stage $t$, and
{let} $u_{-i,t} \in \Delta_{n-1}(\mathcal {A})$ be the empirical
distribution of all consumers' (excluding consumer $i$) actions at
stage $t$.

For a given $n$, it can be seen from (\ref{equa:pricep}) and
(\ref{equa:priceq}) that the prices, and thus the stage payoff in
(\ref{equa:Spayoff1}), are determined by the current global state,
{$\overline{s}_t$,} consumer $i$'s current augmented state,
{$y_{i,t}$,} and current action, {$a_{i,t}$,} as well as the
empirical {distributions, $f_{-i,t}$ and $u_{-i,t}$} of other
consumers' current augmented state and action. Hence, for a certain
function $\pi(\cdot)$, we can write the stage payoff in
(\ref{equa:Spayoff1}) as
\begin{equation}\label{equa:Spayoff}
\pi(y_{i,t},\overline s_t,a_{i,t},f_{-i,t},u_{-i,t}) =
U_t(x_{i,t},s_t,a_{i,t})  - (p_t+w_t) a_{i,t} -   q_t a_{i,t-1} .
\end{equation}

\section{A Continuum Model and Dynamic Oblivious Strategies}\label{sec:DOE}

To study the aggregate behavior of a large number of consumers, we
consider a nonatomic game {involving} a continuum of infinitesimally
small consumers, indexed by $ i \in [0,1]$. We assume that (under state $s_0$) a fraction $\eta_{s_0}$
of the consumers has initial state $x$.
In a nonatomic model,
any single consumer's action has no influence on the aggregate
demand and the prices.
 We consider a class of {strategies (dynamic oblivious
strategies)} in which a consumer's action depends only on the
history of past exogenous states, $h_t=(s_0,\ldots,s_t)$, and her
{own} current state\footnote{Note that a dynamic oblivious strategy
depends only on the consumer's current state,
  instead of her augmented state. As we will see in Section \ref{sec:DOE2},
 in a continuum model, since any single consumer has no influence on the prices,
a best response or equilibrium strategy need not take into account
the action taken at the previous stage.}, i.e., of the form
$$
a_{i,t}=  \overline \nu_{t}( x_{i,t},  h_t).
$$
Suppose that consumer $i$ uses a dynamic oblivious strategy
$\overline \nu=(\overline \nu_0,\ldots,\overline \nu_T)$. Since
there is no idiosyncratic randomness, given a history $h_t$, the
state $x_{i,t}$ of consumer $i$ at stage $t$ depends only on her
initial state $x_{i,0}$. That is, there is a mapping $l_{\overline
\nu,h_t}: \mathcal {X}_0 \to \mathcal {X}_{t}$, such that
$x_{i,t}=l_{\overline \nu,h_t}(x_{i,0})$. Therefore, we can
{specify} the action taken by a dynamic oblivious strategy {in the
alternative form}
 \begin{equation}\label{equa:OE_stra}
a_{i,t}=  \nu_{t}( x_{i,0},  h_t)  \buildrel \Delta \over =
\overline \nu_{t}( l_{\overline \nu,h_t}(x_{i,0}), h_t) .
 \end{equation}
{We refer to} $\nu=(\nu_0,\ldots,\nu_T )$ {as a {\bf dynamic
oblivious strategy}, and let} $\mathfrak{V}$ be the set of all such
strategies.

An alternative {formulation involving strategies that depend on
consumer} expectations on future prices would lead to a Rational
Expectations Equilibrium (REE), an equilibrium concept based on the
rational expectations approach pioneered by \cite{M61}. In our
continuum model, since the only source of stochasticity is from the
exogenous state $s_t$, future prices {under any given strategy
profile,} are completely determined by the history $h_t$. Therefore,
it is {reasonable} to expect that {strategies of the form}
\eqref{equa:OE_stra} will lead to an equilibrium concept that
is identical in outcomes with a REE (cf.\ the discussion in Section
\ref{sec:DOE2}).

Before formally defining a Dynamic Oblivious Equilibrium (DOE), we
{first provide some of the intuition behind the definition.}
 In a continuum model, if all consumers use a common dynamic oblivious strategy $\nu$,
  the aggregate demand and the prices at stage $t$ depend only on the {history} of exogenous states,
$h_t=(s_0,\ldots,s_t)$.  A dynamic oblivious strategy $\nu$ is a DOE
(cf. the formal definition in Section \ref{sec:DOE2}) if it
maximizes every consumer's  {expected total payoff}, under the
sequence of prices that $\nu$ induces. In Section \ref{sec:DOE3}, we
associate a continuum model with a sequence of $n$-consumer models
($n=1,2,\ldots$),
 and
specify the relation between the continuum model and the corresponding $n$-consumer model.

\subsection{The sequence of prices induced by a dynamic oblivious
strategy}\label{sec:DOE1}

Let $h_t=(s_0,\ldots,s_t)$ denote a history up to stage $t$, and
let $\mathcal {H}_t = \mathcal {S}^{t+1}$ denote the set of all
possible such histories. Recall that in a continuum model,
given an initial global state $s_0$, the distribution of consumers'
initial states  is $\eta_{s_0} $. Therefore, under a history $h_t$,
if all consumers use the same dynamic oblivious strategy $\nu$, then
the average demand is
\begin{equation}\label{equa:OE_demand}
\widetilde A_{t \mid \nu, h_{t}}  =  \sum\limits_{x \in \mathcal
{X}_0 } \eta_{s_0} (x) \cdot \nu_{t}(x,h_{t}).
\end{equation}

We now introduce the cost functions in a continuum model. Let
$\widetilde C: \mathbb{R} \times \mathcal {S} \to [0,\infty)$ be a
primary cost function. Let $\widetilde H: \mathbb{R}^2 \times
\mathcal {S}^2 \to [0,\infty)$ be an ancillary cost function at
stage $t \ge 1$, and let $\widetilde H_0: \mathbb{R} \times \mathcal
{S} \to [0,\infty)$ be an ancillary cost function at the initial
stage~$0$.

Given the cost functions in a continuum model, we define the
sequence of prices induced by a dynamic oblivious strategy as
follows:
\begin{equation}\label{equa:OE_price0}
\widetilde p_{t|\nu,h_t} =\widetilde C'(\widetilde
A_{t|\nu,h_t},s_t),\qquad\widetilde q_{0|\nu,h_0}=0, \qquad
\widetilde w_{0|\nu,h_0}= \widetilde H'_0(\widetilde
A_{0|\nu,h_0},s_0),
\end{equation}
and for $t \ge 1$,
\begin{equation}\label{equa:OE_price}
\widetilde q_{t|\nu,h_t}=\dfrac{ \partial \widetilde H\left(
\widetilde A_{t-1|\nu,h_{t-1}}, \widetilde A_{t|\nu,h_{t}},
\overline s_t \right) } { \partial \widetilde
A_{t-1|\nu,h_{t-1}}},
\qquad
\widetilde w_{t|\nu,h_t}=\dfrac{ \partial
\widetilde H\left( \widetilde A_{t-1|\nu,h_{t-1}}, \widetilde
A_{t|\nu,h_{t}}, \overline s_t \right) } { \partial \widetilde
A_{t|\nu,h_{t}}}.
\end{equation}

\subsection{Equilibrium strategies} \label{sec:DOE2}

In this subsection we define the concept of a DOE. Suppose that
all consumers other than $i$ use a dynamic oblivious strategy $\nu$.
In a continuum model, consumer $i$'s action does not affect the
prices. If all consumers except $i$ use a dynamic oblivious strategy
$\nu$,
 consumer $i$'s \textbf{oblivious stage value} (the stage payoff in a continuum model) under a history $h_t$ and an action $a_{i,t}$, is
 \begin{equation}\label{equa:OE_Cpayoffa}
\widetilde \pi_{i,t}(y_{i,t}, h_t, a_{i,t} \mid
 {\nu} )   =
 U_t(x_{i,t},s_t,a_{i,t})  - (\widetilde p_{t|\nu,h_t}+ \widetilde w_{t|\nu,h_t}) a_{i,t} -   \widetilde
 q_{t|\nu,h_t}
 a_{i,t-1},
\end{equation}
where the prices, $\widetilde p_{t|\nu,h_t},$ $\widetilde
w_{t|\nu,h_t}$,  and $\widetilde q_{t|\nu,h_t}$, are defined in
\eqref{equa:OE_price0} and \eqref{equa:OE_price}. Since a single
consumer's action cannot influence $\widetilde q_{t}$, the last term
in \eqref{equa:OE_Cpayoffa} is not affected by the action $a_{i,t}$,
and the decision $a_{i,t}$ at stage $t$ need not take $a_{i,t-1}$
into account, but should take $\widetilde {q}_{t+1}$ into account.

 Consumer $i$'s oblivious stage
value  under a dynamic oblivious strategy $\hat \nu$,
is\footnote{{Recall that} the initial state (the type) of consumer
$i$, $x_{i,0}$, is included in its state $x_{i,t}$, for any $t$.}
 \begin{equation}\label{equa:OE_Cpayoff}
\widetilde \pi_{i,t}(y_{i,t}, h_t  \mid \hat \nu, {\nu} ) \buildrel
\Delta \over = \widetilde \pi_{i,t}(y_{i,t}, h_t, \hat
\nu_t(x_{i,0}, h_t ) \mid
 {\nu} ).
\end{equation}
In particular, we use $\widetilde \pi_{i,t}(y_{i,t},  h_t \mid \nu,
{\nu} ) $ to denote the oblivious stage value of consumer $i$ at
stage $t$, if all consumers use the strategy $\nu$. Given {an
initial global state $s_0$  and an initial state of consumer $i$,
$x_{i,0}$,} her \textbf{oblivious value function} ({total} future
expected payoff function in a continuum model) is
 \begin{equation}\label{equa:OE_Cvalue}
 \widetilde V_{i,{0}} ( x_{i,{0}}, {s_0}  \mid \hat
\nu ,{\nu} ) =
 \E\left\{ {\sum\limits_{\tau={0}}^T { \widetilde
\pi_{i,\tau}(y_{i,\tau},  h_\tau \mid \hat \nu, {\nu} )} } \right\},
 \end{equation}
where the expectation is over the future global states, $\{s_\tau
\}_{\tau=1}^T$.

\begin{definition}\label{Def:OE}
A strategy $\nu$ is a \textbf{Dynamic Oblivious Equilibrium} (DOE)
if
\[
\sup _{\widehat \nu  \in \mathfrak{V} } \;\widetilde V_{i,0}(
x_{i,0}, s_0  \mid \hat \nu, {\nu}) = \widetilde V_{i,0}( x_{i,0},
s_0  \mid  \nu, {\nu}),\qquad\forall x_{i,0} \in \mathcal {X}_{0},
\quad\forall s_0 \in \mathcal {S}.
\]
\end{definition}

A DOE is guaranteed to exist, under suitable assumptions, and this
is known to be the case for our model (under our assumptions),
and even for a more general
model that includes idiosyncratic randomness \citep{BB92}.
 The DOE, as defined above, is
essentially the same concept as the ``dynamic competitive
equilibrium'' studied in \cite{BC10}, which is defined as the
nonatomic equivalent of an MPE, in a continuum model. At a DOE, the
beliefs of all consumers on future prices are consistent with
the equilibrium outcomes. Therefore, a DOE is identical in outcomes
with a Rational Expectations Equilibrium (REE).

In future electricity markets, consumers may form rational
expectations of future prices through an adaptive learning process,
or they may receive price estimates from utilities and/or the
independent system operator through advanced metering
infrastructures. (In Appendix \ref{sec:im}, we {provide} some
discussion {of a possible} implementation of the proposed real-time
pricing mechanism.) {If so, a REE (equivalently, a DOE) will be a
plausible outcome of such a market.} Furthermore, we will show
(Theorem \ref{thm:SO_social}) that under the proposed pricing
mechanism, and under certain convexity assumptions, a DOE is
socially optimal for the continuum model.

For a  dynamic market model with aggregate profit shocks,
\cite{WBV10} introduce a concept of  ``extended oblivious
equilibrium'' at which every firm's strategy depends on its current
state and on the recent history (as opposed to the full history) of
the aggregate shock. The extended oblivious equilibrium is
computationally tractable; however, an equilibrium strategy
 may not be an approximate best response for every firm, even if the
number of firms is large (cf. the error bounds derived in Section
8.3 of \cite{WBV10}).

 Note that the definition of a DOE strategy
requires optimality (attaining the supremum in Definition
\ref{Def:AMP}) only along the equilibrium path \citep{BC10}. Thus, a
DOE is similar in spirit to the ``self-confirming equilibria'' in
\cite{FL93} and the ``subjective equilibria''  in \cite{KL95}, in
which each agent forms correct beliefs about her opponents only
along the equilibrium path.

\subsection{The $n$-consumer model associated with
a continuum model}\label{sec:DOE3}

 We want
the cost functions in a continuum model to approximate the cost
functions in an $n$-consumer model. Since the continuum of consumers
is described by distributions over $[0,1]$, the demand given in
\eqref{equa:OE_demand} can be regarded as the average demand per
consumer. {To capture this correspondence,} we assume the following
relation between the cost functions in a continuum model and their
counterparts in a corresponding $n$-consumer model.

\begin{assumption}\label{A:scale}
For any $n \in \mathbb{N}$, any $s \in \mathcal {S}$, and any
$\overline s$ in $ \mathcal {S}^2$, we have
$$
C^n(A,s) =n \widetilde C \left( \frac{A}{n}, s \right), \;\;\;
 H_0^n(A,s) =n \widetilde H_0 \left( \frac{A}{n}, s \right), \;\;\;
H^n(A,A',\overline s) =n \widetilde H \left( \frac{A}{n},
\frac{A'}{n}, \overline s\right),
$$
where the superscript $n$ is used to indicate that these are the
cost functions associated with an $n$-consumer model.
\end{assumption}

Assumption \ref{A:scale} implies that
$$
(C^n)'(A,s) =\widetilde C'(A/n,s),\qquad (H_0^n)'(A,s) =\widetilde
H_0'(A/n,s),\quad s \in \mathcal {S},
$$
and
$$
\dfrac{\partial H^n(A,A',\overline s)}{\partial A} =\dfrac{\partial
\widetilde H(A/n,A'/n,\overline s)}{\partial
(A/n)},\qquad\dfrac{\partial H^n(A,A',\overline s)}{\partial A'}
=\dfrac{\partial \widetilde H(A/n,A'/n,\overline s)}{\partial
(A'/n)},\;\;\;\forall \overline s \in \mathcal {S}^2,
$$
 i.e., the prices in the continuum model at
the average demand equal the prices \jnt{in} the corresponding
 $n$-consumer model.

\section{Approximation in Large Games} \label{sec:AMP}
In this section, we  consider a sequence of dynamic games, and show that as the
number of consumers increases to infinity,
 a DOE strategy for the corresponding continuum game is asymptotically optimal for every
consumer (i.e., an approximate best response), if the other
consumers follow that same strategy. In the rest of the paper, we
{often} use a superscript $n$ to indicate quantities associated with
an $n$-consumer model.

Suppose that all consumers except $i$ use a dynamic oblivious
strategy $\nu$.  Given a history $h_t$ and an empirical distribution
$f^n_{-i,t}$,
   we use  $v (h_t, f^n_{-i,t}, \nu)$ to denote the empirical
 distribution, {$u^n_{-i,t}$,} of the actions taken by consumers excluding $i$.
In an $n$-consumer model,  suppose that consumer $i$ uses a
{\bf history-dependent strategy} $\kappa^n=\{ \kappa^n_{t} \}_{t=0}^T$ of
the form
\begin{equation}\label{eq:his}
a_{i,t}=\kappa^n_{t}(y_{i,t},h_t,f^n_{-i,t}),
\end{equation}
 while  the other
consumers use a dynamic oblivious strategy $\nu$. {Let}
$\mathfrak{K}_n$ denote the set of all possible history-dependent
strategies $\kappa^n$ for the $n$-consumer model. Note that since
all other consumers use an oblivious strategy $\nu$, $f^n_{-i,t}$ is
completely determined by $\nu$, $f^n_{-i,0}$, and $h_t$.

 The stage payoff received by consumer $i$
at time $t$ is
 \begin{equation}\label{equa:AMP_payoff}
\pi^n_{i,t}( y_{i,t}, h_t, f^n_{-i,t} \mid \kappa^n, {\nu} )= \pi^n
\left(  y_{i,t}, \overline s_t, a_{i,t}, f^n_{-i,t}, v( h_t,
f^n_{-i,t},  {\nu} ) \right),
\end{equation}
where $a_{i,t}=\kappa^n_{t}(y_{i,t},h_t,f^n_{-i,t})$, and the stage
payoff function on the right-hand side is given in
(\ref{equa:Spayoff}).
 Given an initial global state, $s_0$, and consumer $i$'s initial
 state, $x_{i,0}$,
 consumer $i$'s expected payoff under the strategy
$\kappa^n$ is
 \begin{equation}\label{equa:AMP_value}
 V_{i,0} ^n \left(x_{i,0}, s_0 \mid \kappa^n, {\nu}
\right) =  \E\left\{ {\sum\limits_{t=0}^T { \pi_{i,t}^n(  y_{i,t},
h_t, f^n_{-i,t} \mid \kappa^n, {\nu} )} } \right\},
\end{equation}
where the  expectation is over the
 initial distribution  $f^n_{-i,0}$ and over the
future global states, $\{s_t \}_{t=1}^T$. In particular, we use
$V_{i,0} ^n \left(x_{i,0}, s_0 \mid \nu, {\nu} \right)$ to denote
the expected payoff obtained by consumer $i$ if all consumers use
the strategy $\nu$.

\begin{definition}\label{Def:AMP}
A dynamic oblivious strategy $\nu$ has the \textbf{asymptotic Markov
equilibrium (AME)} property \citep{AJ10}, if for any initial global
state $s_0 \in \mathcal {S}$, any initial consumer state $x_{i,0}
\in \mathcal {X}_0 $, and any sequence of history-dependent
strategies $\{\kappa^n\}$, we have
\[
\limsup_{n \to \infty } \left( V_{i,0} ^n \left(x_{i,0}, s_0 \mid
\kappa^n, {\nu}   \right) - V_{i,0} ^n \left( x_{i,0},s_0
 \mid \nu, {\nu} \right) \right) \le 0.
\]
\end{definition}

We will show that every DOE has the
AME property, under the following assumption, which strengthens Assumption \ref{A:price}.

\begin{assumption}\label{A:continuous}
We assume that:
\begin{itemize}
\item [{\rm \ref{A:continuous}.1.}] The {following} four families of functions, {of $A$}, $\{\widetilde
C'({A},s): \; s \in \mathcal{S}\}$, $\{\widetilde H'_0({A}, s): \; s
\in \mathcal{S}\}$, $ \{{\partial \widetilde H(A, A', \overline s)}/{\partial
A}: \; (A',\overline s) \in \mathcal {A} \times \mathcal {S}^2\}$,
and $\{{\partial \widetilde H(A', A, \overline s )}/{\partial A}:
\;(A',\overline s) \in \mathcal {A} \times \mathcal {S}^2\}$, are uniformly
equicontinuous on $[0,\infty)$.\footnote{A sufficient condition for
this assumption to hold is to require a universal bound on the
derivatives of the functions in each family.}

\item [{\rm \ref{A:continuous}.2.}] The marginal costs are bounded from above, i.e.,
 $$
 |\widetilde C'(A,s)| \le P, \qquad  |\widetilde H_0'(A,s)| \le P, \quad \forall (A,s) \in \mathcal {A} \times \mathcal {S},
 $$
 and
 $$
\left|\frac{\partial \widetilde H(A, A', \overline s )}{\partial
A}\right| \le P,\qquad
\left|\frac{\partial \widetilde H(A', A,
\overline s )}{\partial A}\right| \le P,\quad \forall (A',\overline
s) \in \mathcal {A} \times \mathcal {S}^2,
 $$
where $P$ is a positive constant.

\old{\item The transition probability function, $p(x,x',a,s)$, is
continuous in $a$, for any $(x,x',s) \in \mathcal {X}^2 \times
\mathcal {S}. $}

\item [{\rm \ref{A:continuous}.3.}] The utility functions, $\{ U_t(x,s,a) \}_{t=0}^T$, are continuous in $a$ and bounded above, i.e.,
 $$
 U_t(x,s,a) \le Q, \qquad t=0,\ldots,T, \qquad \forall (x,s,a) \in \mathcal {X}_t \times \mathcal {S} \times \mathcal {A},
 $$
where $Q$ is a positive constant.
\end{itemize}
\end{assumption}

Combining with Assumption \ref{A:scale}, Assumption \ref{A:continuous}.1 implies
that  for any \(\varepsilon
>0\), there exists a $\delta>0$ such that for any positive integer
$n$, if $|A-\overline A| \le n \delta$, then
\begin{equation}\label{equa:Acon}
\left|  (C^n)'(A,s) -   (C^n)'(\overline A,s) \right| \le
\varepsilon,\qquad\left|  (H_0^n)'(A,s) -   (H_0^n)'(\overline A,s)
\right| \le \varepsilon,\quad\forall s \in \mathcal {S},
\end{equation}
and for any $(A',\overline s) \in \mathcal {A} \times \mathcal
{S}^2$,
\begin{equation}\label{equa:Acon1}
\left|  \frac{\partial H^n(A, A', \overline s )}{\partial A} -
\frac{\partial H^n(\overline A, A', \overline s )}{\partial
\overline A}\right| \le \varepsilon,\qquad \left|  \frac{\partial
H^n(A', A, \overline s)}{\partial A} - \frac{\partial H^n(
A',\overline A, \overline s )}{\partial \overline A}\right| \le
\varepsilon.
\end{equation}

Note that the boundness of the cost function derivatives implies the
Lipschitz continuity of the cost functions. Combining with
Assumption \ref{A:scale}, for any pair of real numbers $(A,\overline
A)$, and any positive integer $n$, we have
\begin{equation}\label{equa:Alip}
\left|  C^n(A,s) -   C^n(\overline A,s) \right| \le P|A-\overline
A|,\qquad\left| H_0^n(A,s) -   H_0^n(\overline A,s) \right| \le
P|A-\overline A|,\quad\forall s \in \mathcal {S},
\end{equation}
and for any $(A',\overline s) \in \mathcal {A} \times \mathcal
{S}^2$,
\begin{equation}\label{equa:Alip1}
\left|  H^n(A, A', \overline s ) -   H^n(\overline A, A', \overline
s )\right| \le P|A-\overline A|,\qquad  \left| H^n(A', A,  \overline
s ) - H^n(A', \overline A,  \overline s ) \right| \le P|A-\overline
A|.
\end{equation}

\old{Suppose that all consumers other than $i$ use a DOE strategy
$\nu$. We argue (Lemma \ref{L:DP_mu}) that as the number of
consumers increases to infinity, the maximum expected payoff
consumer $i$ can obtain is asymptotically no larger than the optimal
oblivious value, {and that the optimal oblivious value} can be
approximately achieved if she uses the DOE strategy (Lemma
\ref{L:DP_nu}).}

 We argue in the following theorem that a DOE
strategy approximately maximizes a consumer's expected payoff (among
all possible history-dependent strategies) in a dynamic game with a
large but finite number of consumers, if the other consumers also
use that strategy.

\begin{theorem}\label{thm:APE}
Suppose that Assumptions \ref{A:scale}-\ref{A:continuous} hold.
Every DOE has the AME property.
\end{theorem}

Theorem \ref{thm:APE} is proved in Appendix  \ref{sec:A}. Various
approximation properties of nonatomic equilibrium concepts in a
continuum game have been investigated in previous works.
 Sufficient conditions for a
stationary equilibrium (an equilibrium concept for a continuum game
without aggregate uncertainty) to have the AME property are derived
in \cite{AJ10}. For a continuum game with both idiosyncratic and
aggregate uncertainties,  \cite{BC10} shows that as the number
agents increases to infinity, the actions taken in an MPE can be
well approximated by some DOE strategy of the nonatomic limit game.
Note, however, that in a {general} $n$-consumer game, even if all
consumers take an action that is close to the action taken by a DOE
strategy of the nonatomic limit game, the population states and the
prices in the $n$-consumer game can still be very different from
their counterparts in the nonatomic limit game. Therefore, without
further assumptions on the consumers' state transition kernel (e.g., continuous dependence of consumer states on their previous actions),
the
approximation property of a DOE on the action space
 does not necessarily imply the AME property of the DOE.

\section{Asymptotic Social Optimality} \label{sec:social}
In Section \ref{sec:social1}, we define the social welfare
associated with an $n$-consumer model and with a continuum model. In
Section \ref{sec:social2},  we show that for a continuum model, the
social welfare is maximized (over all symmetric dynamic oblivious
strategy profiles) at a DOE, and that for a sequence of $n$-consumer
models, if all consumers use the DOE strategy of the corresponding
continuum model,
  then the
social welfare is asymptotically maximized, as the number of consumers increases to infinity.

\subsection{Social welfare}\label{sec:social1}

In an $n$-consumer model, let {$\v{x}_t=(x_{1,t},\ldots,x_{n,t})$
and $\v{a}_t=(a_{1,t},\ldots,a_{n,t})$ be the vectors of consumer
states and actions, respectively,} at stage $t$. Under the current
{global state $\overline s_t$,} the
 social welfare realized at stage $t$ is
\begin{equation}\label{equa:SO_wel}
 W^n_t(\v{x}_{t},\overline s_t,\v{a}_t) = -C^n(A_t,s_t) -
 H^n(A_{t-1},A_t,\overline s_t) +\sum\limits_{i=1}^n
U_t(x_{i,t},s_t,a_{i,t}),\qquad t=1,\ldots,T,
\end{equation}
and at stage $0$, the social welfare is
\begin{equation}\label{equa:SO_wel0}
 W^n_t(\v{x}_0, s_0,\v{a}_0) = -C^n(A_0,s_0) -
H_0^n(A_0, s_0) +\sum\limits_{i=1}^n U_0(x_{i,0},s_0,a_{i,0}).
\end{equation}

{Because of the symmetry of the problem, the social welfare at stage
$t$ depends on $\v{x}_{t}$ and $\v{a}_t$ only through the empirical distribution of state-action pairs.
In particular, under a symmetric history-dependent
strategy profile $\v{\kappa}^n=(\kappa^n,\ldots,\kappa^n)$ (cf. the
definition of a history-dependent strategy in Eq. \eqref{eq:his}),
we can write the social welfare at time $t$ (with a slight abuse of
notation) as} $W^n_t(f^n_{t},h_t \mid \v{\kappa}^n)$. Given an
initial global state $s_0$ and an initial population state $f^n_0$,
the expected social welfare achieved under a symmetric
history-dependent strategy profile $\v{\kappa}^n$ is given by
\begin{equation}\label{equa:SO_welfarenu}
\mathcal {W}^n_{0}(f^n_0,s_0 \mid \v{\kappa}^n) = W^n_0(f^n_0,s_0
\mid \v{\kappa}^n) +
 \E\left\{ \sum\limits_{t = 1}^T { W^n_t(f^n_t,h_t\mid
\v{\kappa}^n }) \right\},
\end{equation}
where  the expectation is over the future global states
$\{s_t\}_{t=1}^T$. In particular, we use $\mathcal
{W}^n_{0}(f^n_0,s_0 \mid \v{\nu}^n)$ to denote the expected social
welfare achieved by the ``symmetric dynamic oblivious strategy profile'',
$\v{\nu}^n=(\nu,\ldots,\nu)$.

In a continuum model, suppose that all consumers use a {common}
dynamic oblivious strategy $\nu$.
Given an initial global state $s_0$, the expected social welfare is
\begin{equation}\label{equa:SO_CWel}
\widetilde {\mathcal {W}}_0 (s_0 \mid {\nu})=\widetilde W_0( s_0
\mid \nu )+ \E\left\{ \sum\limits_{t=1}^T \widetilde W_t(
h_t \mid \nu) \right\},
\end{equation}
where the expectation is over the future global states,
$\{s_t\}_{t=1}^T$. Here, $\widetilde W_t( h_t \mid \nu)$ is the
stage social welfare under history $h_t$:
\begin{equation}\label{equa:SO_CWelt}
\displaystyle
\begin{array}{l}
 \displaystyle \widetilde W_t( h_t \mid \nu ) =
  -\widetilde C( \widetilde A_{t|\nu,h_t},s_t) -
  \widetilde H(\widetilde A_{t-1|\nu,h_{t-1}},\widetilde A_{t|\nu,h_t},\overline s_t)
  \\[3pt]
\displaystyle \;\;\;\;\;\;\;   \;\;\;\;\;\;\; \;\;\;\;
\;\;\;\;\;\;\;\;\;\; \;\;\;\;\;\;\; \;\;\;\;\;\;\;\;\;\;\;\;\;
\;\;\;\;\; \;\;\;\;\;\; + \sum\limits_{x \in  \mathcal {X}_0 }
\eta_{s_0} (x)
  U_t\left(l_{\nu,h_t}(x),s_t,\nu_t(x, h_t)
  \right),\;\;\;\;\;\;\; \;\;\;t=1,\ldots,T,
\end{array}
\end{equation}
where $l_{\nu,h_t}$ maps a consumer's initial state into her state
at stage $t$, under the history $h_t$ and the dynamic oblivious
strategy $\nu$. The social welfare at stage $0$ is given by
\begin{equation}\label{equa:SO_CWelt0}
 \widetilde W_0(s_0 \mid \nu ) =
  -\widetilde C( \widetilde A_{0|\nu,h_0},s_0) -
  \widetilde H_0(\widetilde A_{0|\nu,h_0}, s_0)
  + \sum\limits_{x \in  {\mathcal {X}}_{0} } \eta_{s_0} (x)
  U_0\left(x,s_0,\nu_0(x,s_0)
  \right).
\end{equation}

 \subsection{Asymptotic social optimality of a DOE}\label{sec:social2}

\old{
\jnt{[IS THIS USEFUL?] Let us start with some perspective on the connection between the continuum and the $n$-consumer model. Let us fix a dynamic oblivious strategy $\nu$ and suppose that it is used by all consumers in the $n$-consumer model, resulting in the social welfare $\mathcal
{W}^n_{0}(f^n_0,s_0 \mid \v{\nu}^n)$. It is not hard to show that as $n\to\infty$, the expected value (over the initial population state  $f^n_0$, which is drawn from the distribution $\eta_{s_0}$) of the social welfare per user, $\mathcal
{W}^n_{0}(f^n_0,s_0 \mid \v{\nu}^n)/n$, converges to the social welfare $\widetilde {\mathcal {W}}_0 (s_0 \mid {\nu})$ under $\nu$ in the continuum model.[IS THIS CORRECT?] The proof (omitted) relies on the laws of large numbers, together with the correspondence
between the two models introduced in Section \ref{sec:DOE3}. Part (b) of Theorem \ref{thm:SO_social} below will be in the same spirit, except that it will make a comparison between with a more general class of strategy profiles (i.e., history-dependent) for the $n$-consumer model and will focus on the issue of optimality.}
}

We now define some notation that will be useful in this subsection.
Since there is no idiosyncratic randomness, given a history $h_t$,
the state of consumer $i$ at stage $t$ depends only on her initial
state $x_{i,0}$, and her actions taken at $\tau=0,\ldots,t-1$.
 At stage $t \ge 1$, the history $h_t$ and the transition function $z_{i,t+1}=r(x_{i,t},a_{i,t},s_{t+1})$ define
a mapping $k_{h_t} : \mathcal {X}_0 \times \mathcal {A}^t  \to
\mathcal {Z}$:
\begin{equation}\label{equa:SO_k}
z_{i,t}=k_{h_t}(x_{i,0},a_{i,0},\ldots,a_{i,t-1}),\;\;\;t=1,\ldots,T.
\end{equation}
Given an initial state $x_{i,0}$,  consumer $i$'s total utility
under a history $h_t$ can be written as a function of her actions
taken at stages $\tau=0,\ldots,t$:
\begin{equation}\label{equa:SO_U}
\overline U_{h_t}(x_{i,0},a_{i,0},\ldots,a_{i,t}) =
U_t(x_{i,0},s_0,a_{i,0})+ \sum\limits_{\tau=1}^t
U_t(x_{i,0},k_{h_\tau}(x_{i,0},a_{i,0},\ldots,a_{i,\tau-1}),s_\tau,a_{i,\tau}).
\end{equation}

Before proving the main result of this section, we introduce a
series of assumptions on the convexity and differentiability of the
cost and the utility functions.

\begin{assumption}\label{A:convex}
We assume the following.
\begin{itemize}

\item [{\rm \ref{A:convex}.1.}] For any $s \in \mathcal {S}$, $\widetilde C(\cdot,s)$ is convex;
for any $\overline s \in \mathcal {S}^2$,  $\widetilde
H(A,A',\overline s)$ is convex in $(A,A')$.

\item [{\rm\ref{A:convex}.2.}] For any $h_T \in \mathcal {H}_T$ and any $x_{i,0} \in \mathcal
{X}_0$, the function defined in (\ref{equa:SO_U}) is concave with
respect to the vector $(a_{i,0},\ldots,a_{i,T})$.

\item [{\rm \ref{A:convex}.3.}] For any $t \ge 1$, any $h_t \in \mathcal {H}_t$, and any $x_{i,0} \in \mathcal
{X}_0$, the function $k_{h_t}$ defined in \eqref{equa:SO_k} is monotonic in
$a_{i,\tau}$, for $\tau=0,\ldots,t-1$; further, its left and right
derivatives with respect to  $a_{i,\tau}$ exist, for
$\tau=0,\ldots,t-1$.

\item [{\rm \ref{A:convex}.4.}] For $t \ge 1$, and for any $(x,s,a) \in \mathcal {X}_0 \times \mathcal {S} \times \mathcal {A} $,
 the left and right derivatives of the utility function $U_t(x,z,s,a)$ in $z$ exist.
\end{itemize}
\end{assumption}

Assumption \ref{A:convex}.1 is standard. If the utility function is
concave in $a$, Assumption \ref{A:convex}.2 requires that the
transition function $k_{h_t}$ preserves concavity (a linear function would be an example).
 Note that Assumptions
\ref{A:convex}.1 and \ref{A:convex}.2  guarantee that in both models
(a dynamic game with a finite number of consumers, and the
corresponding continuum game), the expected social welfare (consumer
$i$'s expected payoff) is concave in the vector of actions taken by all
consumers (respectively, by consumer $i$).
 Assumptions \ref{A:convex}.3 and
\ref{A:convex}.4 ensure the existence of left and right derivatives
of the expected social welfare given in (\ref{equa:SO_CWel}), with
respect to the actions taken by consumers. An example where
Assumptions \ref{A:convex}.2-\ref{A:convex}.4  hold is given next.

\begin{example}{\rm
Consider appliances such as Plug-in Hybrid Electric
Vehicles (PHEVs), dish washers, or clothes washers. For such
appliances, a customer usually only cares  whether a task is completed
before a certain time.

 Given an initial state (type) of consumer
$i$, $x_{i,0}$, let $D(x_{i,0})$ and $T(x_{i,0})$ indicate her total
desired demand and the stage {by which} the task has to be
completed, respectively. Under a given history $h_t$, the total
utility {accumulated} by consumer $i$ until time $t$ is assumed to
be of the form
$$
\overline U_{h_t}(x_{i,0},a_{i,0},\ldots,a_{i,t}) = Z\left(x_{i,0},
\min\left\{ D(x_{i,0}), \sum\limits_{\tau=0}^{\min\{T(x_{i,0}),t\}}
a_{i,\tau} \right\} \right),
$$
for some function $Z$. If for every $x_{i,0} \in \mathcal {X}_0$,
$Z(x_{i,0},\,\cdot\,)$ is nondecreasing and concave, then Assumption
\ref{A:convex}.2 holds. At stage $t=0$, we have
$$
U_{0}(x_{i,0},s_0,a_{i,0})=Z\left( x_{i,0}, \min\left\{ D(x_{i,0}),
a_{i,0} \right\} \right).
$$
For $t=1,\ldots, T(x_{i,0})$, we {let} $z_{i,t}=
\sum\nolimits_{\tau=0}^{t-1} a_{i,\tau}$, and
$$
U_t(x_{i,0},z_{i,t},s_t,a_{i,t})=Z\left( x_{i,0}, \min\left\{
D(x_{i,0}), a_{i,t}+ z_{i,t} \right\} \right)-Z\left( x_{i,0},
\min\left\{ D(x_{i,0}), z_{i,t} \right\} \right) .
$$
For $t \ge T(x_{i,0})+1$, we let $z_{i,t}= D(x_{i,0})$, and let $
U_t(x_{i,t},s_t,a_{i,t}) $ be identically zero. Suppose that for
every $x_{i,0} \in \mathcal {X}_0$, the right and left derivatives
of $Z(x_{i,0},\,\cdot\,)$ exist. Then, Assumptions \ref{A:convex}.3 and
\ref{A:convex}.4 hold. }$\;\;\;\square$
\end{example}

\begin{theorem}\label{thm:SO_social}
Suppose that Assumptions \ref{A:scale}-\ref{A:convex} hold. Let
$\nu$ be a DOE of the continuum game. Then, the following hold.
\begin{itemize}
\item[(a)] In the continuum game, the social welfare is maximized
{(over all symmetric dynamic oblivious strategy profiles)} at the
DOE, i.e., \footnote{{Note that we are only comparing the social
welfare under different {\bf symmetric} dynamic oblivious strategy
profiles, where all consumers are using the same dynamic oblivious
strategy ($\nu$ or $\vartheta$). This is no loss of generality
because} under Assumption \ref{A:convex},
  the social welfare in a continuum game is {a concave function of the collection of consumer actions taken under the} different histories.
  Hence, {it can be shown that}
 the optimal social welfare can be achieved by a symmetric dynamic oblivious strategy
profile.}
$$
  \widetilde {\mathcal W}_0(s_0 \mid {\nu})  = \sup\nolimits_{ \vartheta \in \mathfrak{V}}
   \widetilde {\mathcal W}_0(s_0 \mid \vartheta)  ,
 \;\;\; \forall s_0 \in \mathcal {S},
$$
where $\mathfrak{V}$ is the set of all dynamic oblivious strategies.

\item[(b)] For a sequence of $n$-consumer games, the symmetric DOE strategy profile, $\v{\nu}^n=(\nu,\ldots,\nu)$, approximately maximizes the
 expected social welfare, as the number of consumers increases to infinity. That is,
 for any
initial global state $s_0$, and any sequence of symmetric
history-dependent strategy profiles $\{\v{\kappa}^n\}$, we
have\footnote{Under Assumption \ref{A:convex},  the social welfare
in an $n$-consumer game is {a concave function of the collection of
consumer  actions taken under the} different histories. Therefore,
$\sup_{{\kappa}^n \in \mathfrak{K}_n} \mathcal {W}^n_{0}(f^n_0,s_0
\mid \v{\kappa}^n) $ is also the maximum social welfare that can be
achieved by a (possibly non-symmetric) history-dependent strategy
profile. }
\[
 \limsup_{n \to \infty } \E\left\{
\dfrac{   \mathcal {W}^n_{0}(f^n_0,s_0 \mid \v{\kappa}^n)-\mathcal
{W}^n_{0}(f^n_0,s_0 \mid \v{\nu}^n)}{n}   \right\} \le 0,
\]
where the expectation is over the initial population state, $f^n_0$.
\end{itemize}

\end{theorem}

The proof of Theorem \ref{thm:SO_social} is given in Appendix
\ref{sec:C}.

\section{Conclusion and Future Directions}\label{sec:con}

In an electric power system, load swings may result in significant
ancillary cost to suppliers. {Motivated} by the observation that
  marginal cost pricing  may not achieve social optimality in
electricity markets, we proposed a new dynamic pricing mechanism
that takes into account the externality conferred by a consumer's
action on future ancillary cost. Besides proposing a suitable
game-theoretic model that incorporates the cost of load fluctuations
and
 a particular pricing mechanism for electricity markets, a main contribution of this paper
was to show that the proposed pricing mechanism achieves social
optimality in a dynamic nonatomic game, and approximate social
optimality for the case of finitely many consumers, under certain
convexity assumptions.

To compare the proposed pricing mechanism with marginal cost
pricing, we presented a numerical example in which the demand
increases sharply at the last stage. In this example,  the proposed
pricing mechanism creates a stronger incentive for consumers to
shift their peak load than marginal cost pricing, through an
additional negative price charged on {off-peak consumer demand.} As
a result, compared with marginal cost pricing, the proposed pricing
mechanism achieves a higher social welfare, and at the same time,
reduces the peak load, and therefore has the potential to reduce the
need for long-term investments in peaking plants.

We believe that the constructed dynamic game-theoretic model, the
proposed pricing mechanism, and more importantly, the insights
provided by this work, can be applied to a more general class of
markets with friction.  As an extension and future work, one can
potentially develop and use variations of our framework to a market
of a perishable product/service where demand fluctuations incur
significant cost to suppliers. {Examples include} data centers
implementing cloud services that suffer from the switching costs to
toggle a server into and out of a power-saving mode \citep{LW11},
and large organizations such as hospitals that use on-call staff to
meet unexpected demand.

\bibliographystyle{nonumber}

\ACKNOWLEDGMENT{The authors are grateful to Prof.\ Michael Caramanis
for discussions and several comments on a draft of this paper, and to Prof. Ramesh Johari for discussions
and pointers to the literature. This
research was supported in part by the National Science Foundation
under grant CMMI-0856063 and by a Graduate Fellowship from Shell.}

\newpage


\begin{appendices}
\textbf{\large{Electronic Companion}}
\section{Approximation of the supplier cost}\label{sec:A0}

In this appendix, we
show via simulation that at least in some cases,
the supplier cost (including the cost of ancillary service) can be captured by
 a simplified
cost function of the form in \eqref{equa:cost}. We consider a
$(T+1)$-stage dynamic model with two energy resources, a primary
energy resource and an ancillary energy resource. It is assumed that
the forecast demand is met by the primary energy resource (e.g.,
coal-fired or nuclear power generators), and that at stage
$t=1,\ldots,T$,  the deviations from the forecast demand,
$\{w_t\}_{t=1}^T$, are independent random variables uniformly
distributed on $[-\omega,\omega]$. At the initial stage $0$, we
assume that the forecast error is zero, i.e., $w_0=0$.

At stage $t$, let $b_t$ denote the difference between the actual
output of the primary energy resource and the forecast demand, and
let $d_t$ denote the output of the ancillary energy resource (e.g.,
oil/gas combustion turbines). For simplicity, we {will} assume that
the cost of {a positive primary energy resource (respectively,
ancillary energy resource) is $b^2_t$ (respectively, $10 d^2_t$).}

Let $r_b$ be the ramping rate of the primary energy resource, and
$r_d$ be the ramping rate of the ancillary energy resource. At the
initial stage $0$, we assume that $b_0=w_0=0$, and $d_0=0$.
 At stage $t \ge 1$, if $w_t<0$, then $d_t=0$, and we assume that $b_t=0$, that is,
 the system operator maintains a high level of (potential) output in order to be able to deal with a possible unexpected
 demand surge  in the future; if
$w_t>0$, we assume that $b_t= \min \{w_t, b_{t-1}+r_b\}$, where
$b_{t-1}+r_b$ is the maximum possible output of the primary energy
resource at stage $t$, and that $d_t = \min\{ w_t-b_t, d_{t-1}+r_d
\}$. The total supplier cost (excluding the cost to meet the
forecast demand) is
\begin{equation}\label{equa:C}
C=\sum_{t=1}^T \left( b^2_t+ 10 d^2_t \right).
\end{equation}
For notational convenience, we let $(\cdot)^+=\max\{ \cdot,0\}$. We
use the following function to approximate the supplier cost:
\begin{equation}\label{equa:AC}
\widetilde C=\sum_{t=1}^T \left( \widetilde b^2_t+ 10 \widetilde d^2_t\right),
\end{equation}
where $\widetilde d_t=\min \left\{  r_d,
\big(0,w_t-(w_{t-1})^+-r_b \big)^+ \right\}$, and $\widetilde
b_t=(w_t-\widetilde d_t)^+$.

The function in  \eqref{equa:AC} well approximates the supplier cost
in  \eqref{equa:C}, if for an unexpected demand surge at stage $t$,
the system load at the previous stage, $w_{t-1}$, is met by the
primary energy resource, and load shedding rarely occurs (so that
$(w_t)^+$ typically equals $b_t+d_t$). Note that in \eqref{equa:AC},
for each stage  $t$, the approximated cost depends only on $w_{t-1}$
and $w_t$. Therefore, the approximated cost in \eqref{equa:AC} can
be written as
\begin{equation}\label{equa:AC1}
\widetilde C=\sum_{t=1}^T  \left( ((w_t)^+)^2 + H(w_{t-1},w_t)\right),
\end{equation}
where $H(w_{t-1},w_t) =\left(\widetilde b^2_t+ 10 \widetilde d^2_t-((w_t)^+)^2\right)^+$.

 \begin{figure}\label{Fig:app}
    \includegraphics[width=10cm]{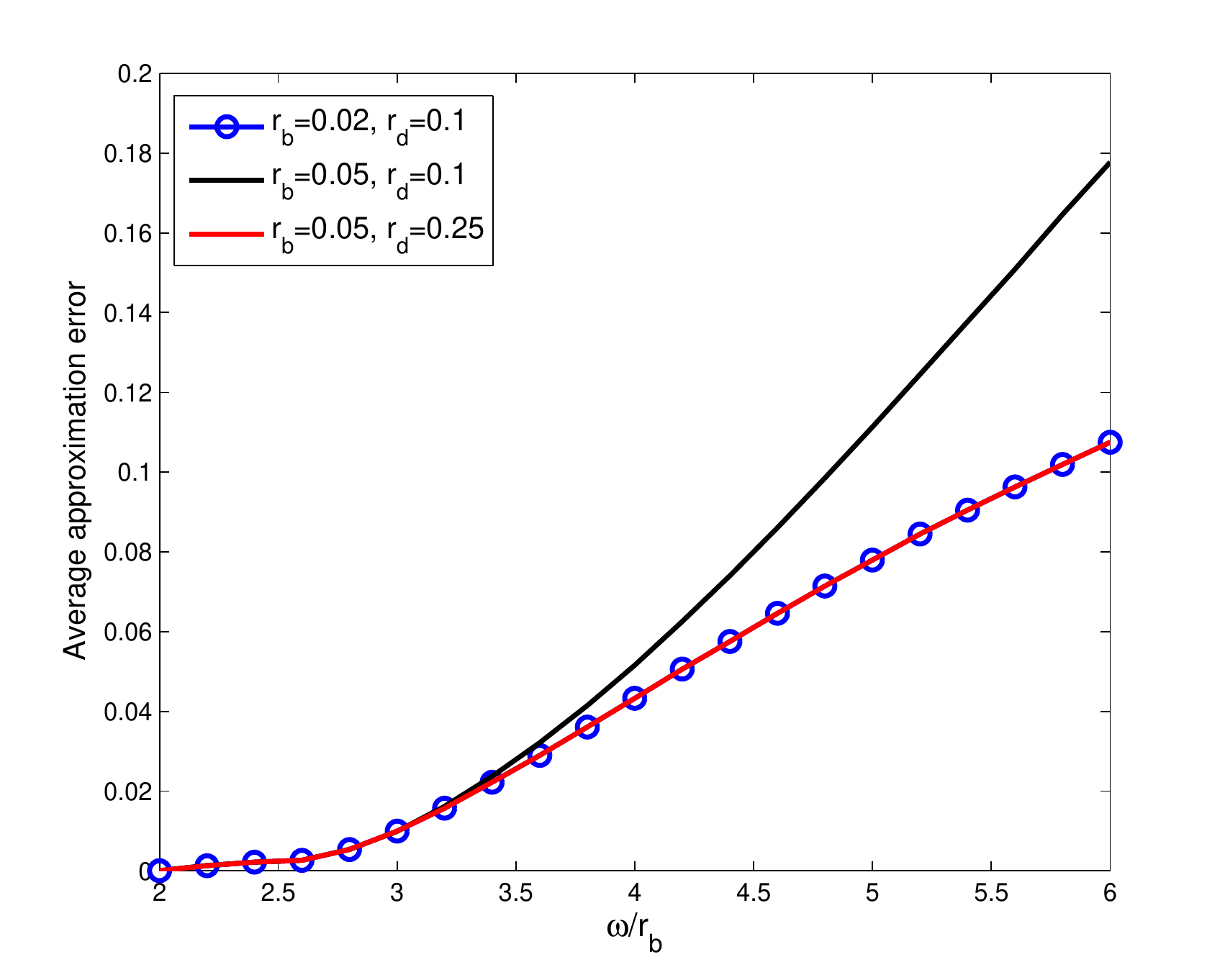}
    \centering
    \caption{A simulation experiment with $T=24$, and $500,000$ trajectories for each $\omega/r_b$ on the horizontal axis. The approximation
    error is defined by
    $|C- \widetilde C|/C$. The average approximation
    error (vertical axis) is the mean value of the approximation
    errors of the $500,000$ trajectories.}
  \end{figure}

For different values of the parameters, $r_b$, $r_d$, and $\omega$,
we evaluate the performance of the approximation via simulation.
Fig.\ 2 depicts some numerical results of a simulation experiment
and we can make the following
 observations:
\begin{enumerate}
\item The main source of approximation error is from the following scenario: at stage $t-1$,
 the deviation in demand $w_{t-1}$ is nonpositive,
 $w_t>r_b$, and  $w_{t+1}>2r_b$.
 In this scenario, the output of the primary energy source at stage $t$ is $r_b$, which is less than $w_{t}$.
  When $\omega/r_b \le 2$, this scenario never occurs and we observe from Fig.\ 2 that the
  approximation error is close to zero, regardless of the value of $r_d$.

\item Comparing the black curve with the red curve in Fig.\ 2, we observe that when $\omega/r_b>3$ (when $r_b=0.05$ and $\omega > 0.15$),
the approximation error for the case where $r_d=0.1$ is larger than
that for the case where $r_d=0.25$. This is because for the case with $r_d=0.1$,
as $\omega/r_b$ increases from $3$ to $6$ (as $\omega$ increases from $0.15$ to $0.3$),
the probability of load shedding increases, which deteriorates the performance of the approximation.

\item Finally, and perhaps most importantly, when the ramping rate of the ancillary energy resource is high enough to prevent any load shedding,
the approximation error is an increasing function of the single
parameter $\omega/r_b$ (e.g., the blue curve with circle markers for
$r_b=0.02, r_d=0.1$ and the red curve for $r_b=0.05, r_d=0.25$ merge
together in Fig.\ 2); in this case, we observe from Fig.\ 2 that the
approximation error is less than $10 \%$ for a wide range of
parameter values.

\end{enumerate}

\section{Implementation of the proposed pricing mechanism}\label{sec:im}

To implement a dynamic real-time pricing mechanism, all consumers
should be exposed to time-varying prices associated with ex ante
estimates of generation costs that reflect system operating
conditions (p. 81 of \cite{BJR02}) so that they can adjust
 their demand in accordance to real-time prices as well as ex ante price estimates.
  The mechanism proposed in this paper is not an exception.
The ex ante estimates of real-time prices can be developed by
evaluating statistical relationships between historical real-time
prices and various factors such as load forecast, weather
predictions, and expected supply/demand balances \citep{BJR02}.

We now provide a brief discussion of the details of a possible
implementation of the proposed pricing mechanism:
\begin{description}{\itemsep=0pt}

\item[Ex ante price estimates.]  Suppose that the exogenous state $s_t$ is realized at the
beginning of each stage $t$; for every possible realization of {the
trajectory (scenario) of} future exogenous states
$\{s_\tau\}_{\tau=t+1}^{t+T}$, consumers receive corresponding
 price estimates  $\{\hat p_\tau\}_{\tau=t}^{t+T}$,
$\{\hat w_\tau\}_{\tau=t}^{t+T}$, and $\{\hat
q_\tau\}_{\tau=t}^{t+T}$,
 from utilities and/or the independent
system operator. {The
consumers also know or receive the probabilities of the different
trajectories.} With the received price estimates (associated with
possible {trajectories} of future exogenous states) and preset
utility functions, each consumer's infrastructure solves a dynamic
programming problem to maximize her expected payoff  over the
horizon from $t$ to $t+T$. (The state at time $\tau$ in this dynamic
program is comprised of $y_{i,\tau}$, and the history
$(s_t,s_{t+1},\ldots,s_{\tau}$.) The dimension of this state space
grows with the time horizon $T$ (because of the exponentially
increasing number of histories). Unfortunately, this is unavoidable
for models of this type, and might require some further
approximations, e.g., in the spirit of \cite{WBV10}.

\item[Ex post prices.] At each stage $t$,  after the realization of the system
demands
$A_{t-1}$ and $A_t$, consumers pay ex post prices $(p_t,w_t,q_t)$
that are determined according to Eqs. \eqref{equa:pricep} and
\eqref{equa:priceq}.

\item[Equilibrium.]
 In a
market with a large number of price-taking consumers, it is possible
to make ex ante price estimates (contingent on the realized
trajectories) that are close to ex post prices. If every consumer
maximizes her own payoff in response to these pretty accurate price
estimates, the resulting outcome {should} be close to that resulting
from a Rational Expectations Equilibrium (REE). The results derived
in this paper show that the expected social welfare can be
approximately maximized, under the proposed mechanism.
\end{description}

We emphasize here that there remain several challenging implementation issues, e.g., the accuracy
of future price estimates and the uncertainty of consumer response to ex ante
price estimates. {For example,} \cite{RDM10} show that if consumers
act myopically to highly inaccurate price estimates, real-time
pricing may result in extreme price volatility. However, we note
that these challenges are {generic to} almost all kinds of real-time
pricing mechanisms.

\section{Proof of Theorem \ref{thm:APE}}\label{sec:A}

We consider a sequence of $n$-consumer models where $n-1$ consumers
(all except for consumer $i$) use a DOE strategy $\nu$.
 As the number of
consumers increases to infinity, the randomness of consumer initial
states averages out. Thus, in Step 1 we show that the aggregate
demand (in an $n$-consumer model) at a history $h_t$ is close to
$n\widetilde A_{t|\nu,h_t}$ (defined in Eq. \eqref{equa:OE_demand}), with high probability. As a
consequence, we show in Step 2 that as $n\to\infty$, consumer $i$'s
expected payoff associated with any sequence of actions can be
approximated by her oblivious value defined in
\eqref{equa:OE_Cvalue}. Since the DOE strategy $\nu$ maximizes
consumer $i$'s oblivious value among all possible strategies, we
argue in Step 3 that as $n\to\infty$, the maximum expected payoff
consumer $i$ can obtain is asymptotically no larger than the optimal
oblivious value. In Step 4, we show that consumer $i$'s optimal
oblivious value can be approximately achieved if she uses the DOE
strategy $\nu$. We finally conclude with the AME property of the DOE
strategy $\nu$ (the result in Theorem \ref{thm:APE}).

In what follows, we will be using the uniform metric over the set of
probability distributions on the finite set $\mathcal {X}_0$.
Specifically, if $f$ and $f'$ are two distributions on ${\mathcal
{X}}_{0}$, we let
\begin{equation}\label{equa:SO_dis}
d(f,f')  \buildrel \Delta \over
=\left\|f-f'\right\|_{\infty}= \max_{x \in  {\mathcal {X}}_{0}}
\left| f(x)-f' (x) \right|.
\end{equation}

\textbf{Step 1:} {\it With high probability, the aggregate demand
{under} a history $h_t$ is close to $n\widetilde A_{t|\nu,h_t}$.}

Given an initial distribution $f^n_{-i,0}$, {and if} all consumers
(excluding $i$) use a dynamic oblivious strategy $\nu$, we write
their aggregate demand at a history $h_t$ as
$$
A^n_{-i,t}=(n-1)\sum\limits_{x \in \mathcal {X}_0} f^n_{-i,0}(x)
\nu_t(x, h_t).
$$
Recall that (cf. \eqref{equa:OE_demand})
$$ \widetilde A_{t \mid \nu, h_{t}}  =  \sum\limits_{x
\in \mathcal {X}_0 } \eta_{s_0} (x) \cdot \nu_{t}(x,h_{t}).
$$
We observe that if $d(f^n_{-i,0},\eta_{s_0}) \le \delta/(XB)$, then
at any history $h_t$ we have
\begin{equation}\label{equa:AMP_A}
\left| A^n_{-i,t} - (n-1) \widetilde A_{t|\nu,h_t}  \right| \le
\delta (n-1),
\end{equation}
with probability at least $1-O(e^{-n})$. {More precisely,} since
{the} consumers' initial states are independently drawn according to
$\eta_{s_0}$, {Hoeffding's inequality (\cite{H63}) yields,}
\begin{equation}\label{equa:AMP_Hoe}
\mathbb{P}\left( d(F_{s_0}^{n-1},\eta_{s_0}) \ge \delta/(XB) \right)
\le 2 X \exp \left\{-2(n-1)\delta^2/(X^2B^2)\right\},\;\;\forall s_0
\in \mathcal {S},\;\;\forall \delta>0,\;\;\forall n \in
\mathbb{N}^+,
\end{equation}
where {$X$ is the cardinality of the set $\mathcal {X}_0$ and}
$F_{s_0}^{n-1}$ is an $X$-dimensional  random vector denoting the
distribution of the initial states of the $n-1$ consumers (excluding
$i$).

\bigskip

\textbf{Step 2:} {\it Under a given history $h_T$, consumer $i$'s
expected payoff can be approximated by a corresponding oblivious
value, defined in \eqref{equa:L_ob}.}

In an $n$-consumer model, suppose that all consumers other than $i$
use a dynamic oblivious strategy $\nu$. Given a complete history
$h_T=(s_0,\ldots,s_T)$, and consumer $i$'s initial state $x_{i,0}$,
{we define} her expected payoff under a history-dependent strategy
$\kappa^n$ by
$$
V^n_{i,0}(x_{i,0},h_T \mid \kappa^n, \nu) = \E \left\{
V^n_{i,0}(x_{i,0},h_T,f^n_{-i,0} \mid \kappa^n, \nu) \right\} ,
$$
where the expectation is over the initial distribution,
$f^n_{-i,0}$, and $V^n_{i,0}(x_{i,0},h_T,f^n_{-i,0} \mid \kappa^n,
\nu) $ is consumer $i$'s payoff under the given initial distribution
$f^n_{-i,0}$,
\begin{equation}\label{equa:L_re}
V^n_{i,0}(x_{i,0},h_T,f^n_{-i,0} \mid \kappa^n,
\nu)={\sum\limits_{t=0}^T {
\pi_{i,t}^n(  y_{i,t}, h_t, f^n_{-i,t} \mid \kappa^n, {\nu} )} },
\end{equation}
and where the stage payoff function, $\pi_{i,t}^n(\cdot)$, {has been
defined} in \eqref{equa:AMP_payoff}. Note that given $f^n_{-i,0}$,
and since all consumers other than
 $i$ use a dynamic oblivious strategy, the distribution of their
 augmented states, $f^n_{-i,t}$, is {completely} determined by the history $h_t$.
Therefore, given $f^n_{-i,0}$, consumer $i$'s history-dependent
strategy $\kappa^n$ is equivalent to a dynamic oblivious strategy:
the action it takes at stage $t$ depends only on $x_{i,0}$ and
$h_t$. We can therefore define an oblivious strategy $\widetilde
\nu^n(\kappa^n, f^n_{-i,0})$ such that
$$
\widetilde \nu_t(\kappa^n, f^n_{-i,0}) (x_{i,0},h_t) =
\kappa^n_{t}(y_{i,t}, h_t, f^n_{-i,t}),
$$
where $f^n_{-i,t}$ is the distribution of the $n-1$ consumers'
augmented states under the history $h_t$, induced from the initial
distribution $f^n_{-i,0}$ by the symmetric oblivious strategy
profile $(\nu,\ldots,\nu)$, and $y_{i,t}$ is consumer $i$'s
 augmented state under the history $h_t$, induced from her initial state $x_{i,0}$ by the strategy $\kappa^n$.

In the corresponding continuum model, suppose that all consumers
other than $i$ use a dynamic oblivious strategy
 $\nu$. For a given complete history $h_T$, we
 define consumer $i$'s oblivious value under
 an initial distribution $f^n_{-i,0}$, her initial state $x_{i,0}$,
 and the history-dependent strategy $\kappa^n$:
\begin{equation}\label{equa:L_ob1}
 \widetilde V_{i,0} ( x_{i,0}, h_T, f^n_{-i,0} \mid \kappa^n, \nu ) =  {\sum\limits_{t=0}^T { \widetilde
\pi_{i,t}(y_{i,t}, h_t \mid  \widetilde \nu(\kappa^n, f^n_{-i,0}) ,
{\nu} )} } ,
\end{equation}
where the oblivious stage value function $\widetilde
\pi_{i,t}(\cdot)$ is given in \eqref{equa:OE_Cpayoff}. We define the
expected oblivious value for consumer $i$ under the
history-dependent strategy $\kappa^n$, as\footnote{This is actually
the oblivious value achieved by a mixed strategy under the complete
history $h_T$. In the continuum model, under a history $h_t$,  the
mixed strategy takes an action $\widetilde \nu_t(\kappa^n,
f^n_{-i,0}) (x_{i,0},h_t)$, if the distribution of the $n-1$
consumers' (excluding $i$'s) initial states in the corresponding
$n$-consumer model is realized as $f^n_{-i,0}$.}
\begin{equation}\label{equa:L_ob}
 \widetilde V_{i,0} ( x_{i,0},h_T \mid \kappa^n, \nu ) =  \E \left\{
 \widetilde V_{i,0} ( x_{i,0}, h_T, f^n_{-i,0} \mid \kappa^n, \nu )
 \right\},
\end{equation}
where the expectation is over the initial distribution,
$f^n_{-i,0}$. For any \(\varepsilon  > 0\), in this step we aim to
show that there exists a positive integer $N$ such that for any
sequence of history-dependent strategies $\{\kappa^n\}$,
\begin{equation}\label{equa:L_21}
\left| \widetilde V_{i,0} ( x_{i,0},h_T  \mid \kappa^n, \nu ) -
V^n_{i,0}(x_{i,0},h_T \mid \kappa^n, \nu) \right| \le \varepsilon,
\;\;\;\forall n \ge N, \;\;\;\forall h_T \in \mathcal {H}_T,\;\;\;
\forall x_{i,0} \in \mathcal {X}_0.
\end{equation}

 For a given $s_0$, let $\mathfrak{F}^{n-1}_{s_0}(\delta)$
be the set of $f^n_{-i,0}$ such that $d(f^n_{-i,0},\eta_{s_0} ) \le
\delta$. To verify \eqref{equa:L_21}, we first argue that for any
\(\varepsilon  > 0\), there exists an positive integer $N_1$ and
some $\delta>0$ such that for any $f^n_{-i,0} \in
\mathfrak{F}^{n-1}_{s_0}(\delta/(XB))$ and any $n \ge N_1$,
\begin{equation}\label{equa:L_211}
\left| \widetilde V_{i,0} ( x_{i,0}, h_T, f^n_{-i,0} \mid \kappa^n,
\nu ) - V^n_{i,0}(x_{i,0}, h_T, f^n_{-i,0} \mid \kappa^n, \nu)
\right| \le \varepsilon/2,  \;\;\;\forall h_T \in \mathcal
{H}_T,\;\;\; \forall x_{i,0} \in \mathcal {X}_0.
\end{equation}

Under the uniform equicontinuity assumption for the derivatives of the cost
functions (see Eqs. \eqref{equa:Acon} and \eqref{equa:Acon1}), we
know that a small deviation of the aggregate demand from $\widetilde
A_{t|\nu,h_t}$ will result in  prices that are only slightly
different from the prices in the continuum model. We also note that
consumer $i$ cannot take an action larger than $B$, and her payoff
is influenced by other consumers only through the prices. For any
$\varepsilon>0$, we can find some $\delta>0$ and a positive integer
$N_1$ such that for any given
 $(x_{i,0}, h_t)$, if $f^n_{-i,0} \in
\mathfrak{F}^{n-1}_{s_0}(\delta/(XB))$, then the inequality in
\eqref{equa:AMP_A} holds for any history $h_\tau$, which implies
that
\begin{equation}\label{equa:AMP_Lmu11}
\left|  \widetilde \pi_{i,t}(y_{i,t}, h_t, \kappa^n_t(y_{i,t}, h_t,
f^n_{-i,t}) \mid  {\nu} )-  \pi_{i,t}^n( y_{i,t}, h_t, f^n_{-i,t}
\mid \kappa^n, {\nu} )   \right| \le \varepsilon/(2T+2),\;\;\forall
n \ge N_1,\;\;\forall h_t,
\end{equation}
i.e., consumer $i$'s stage payoff (under the action
$\kappa^n_t(y_{i,t}, h_t, f^n_{-i,t})$) in the $n$-consumer model is
close to her oblivious stage value (under the same action
$\kappa^n_t(y_{i,t}, h_t, f^n_{-i,t})$) in the continuum model, if
the initial distribution in the $n$-consumer model, $f^n_{-i,0}$, is
close to its expectation. The result in \eqref{equa:L_211} follows
from Eq. \eqref{equa:AMP_Lmu11} and the definitions in
\eqref{equa:L_re} and \eqref{equa:L_ob1}.
  Note that
$Q+2BP$ is an upper bound on the stage payoff that consumer $i$
could obtain, and $-2BP$ is a lower bound on consumer $i$'s stage
payoff, under Assumption \ref{A:continuous}. The desired result in
\eqref{equa:L_21} follows from \eqref{equa:L_211}, and the fact that
the probability that $f^n_{-i,0} \notin
\mathfrak{F}^{n-1}_{s_0}(\delta/(XB))$ decays exponentially with $n$
(cf. Eq. \eqref{equa:AMP_Hoe}).

\bigskip

\textbf{Step 3:} {\it The maximum expected payoff consumer $i$ can
obtain is asymptotically no larger than the optimal oblivious
value.}

In this step, we consider the case where all consumers in an
$n$-consumer model except for $i$ use a DOE strategy $\nu$, and
argue that for any sequence of history-dependent strategies
$\{\kappa^n\}$,
\begin{equation}\label{equa:dynamic_mu}
 \limsup \limits_{n \to \infty }  \left( V_{i,0} ^n
\left(x_{i,0},s_0 \mid \kappa^n, {\nu}   \right) - \widetilde
V_{i,0} ( x_{i,0}, s_0 \mid {\nu}, {\nu}) \right) \le 0, \;\;\;
\forall s_0 \in \mathcal {S}, \;\;\; \forall x_{i,0} \in \mathcal
{X}_0,
\end{equation}
where consumer $i$'s expected payoff, $V_{i,0} ^n \left(x,s \mid
\kappa^n, {\nu}   \right)$, is given in (\ref{equa:AMP_value}), and
$\widetilde V_{i,0} ( x, s \mid {\nu}, {\nu})$ is the oblivious
value function in (\ref{equa:OE_Cvalue}). We first observe that
\begin{equation}\label{equa:L_24}
  V_{i,0} ^n
\left(x_{i,0},s_0 \mid \kappa^n, {\nu}   \right) = \sum\limits_{h_T
\in\mathcal {H}_{T}(s_0) } \mathbb{P}(h_T\mid s_0) \cdot V_{i,0} ^n
\left(x_{i,0},h_T \mid \kappa^n, {\nu} \right),
\end{equation}
where $\mathcal {H}_{T}(s_0)$ is the set of complete histories
commencing at state $s_0$, and $ \mathbb{P}(h_T\mid s_0)$ is the
probability that the history $h_T$ is realized, conditional on the
initial global state being $s_0$. We define
\begin{equation}\label{equa:L_224}
\widetilde V_{i,0} ( x_{i,0}, s_0 \mid \kappa^n, {\nu} )  =
\sum\limits_{h_T \in\mathcal {H}_{T}(s_0) } \mathbb{P}(h_T\mid s_0)
\cdot \widetilde V_{i,0} ( x_{i,0},h_T  \mid \kappa^n, \nu ).
\end{equation}
Note that if $\kappa^n$ happens to be a dynamic oblivious strategy,
this definition is consistent with the definition of oblivious value
function in \eqref{equa:OE_Cvalue}.

For any \(\varepsilon  > 0\), let $N$ be the integer defined in Eq.
\eqref{equa:L_21}; for any sequence of history-dependent strategies
$\{\kappa^n\}$, we argue that
\begin{equation}\label{equa:L_255}
\begin{array}{l}
\displaystyle \widetilde V_{i,0} ( x_{i,0}, s_0 \mid {\nu},{\nu})
\ge \sum\nolimits_{h_T \in \mathcal {H}_{T}(s_0)} \mathbb{P}(h_T\mid
s_0) \cdot \widetilde V_{i,0}
\left(x_{i,0},h_T \mid \kappa^n, {\nu} \right)\\[5pt]
\displaystyle \;\;\;\;\;\;\;\;\;\;\;\;\;\;\;\;\;\;\;\;\;\;\;\;\;
\;\;\;\;\;\;\;\; \;\;\; \;\ge \sum\nolimits_{h_T \in \mathcal
{H}_{T}(s_0)} \mathbb{P}(h_T\mid s_0) \cdot (V^n_{i,0}(x_{i,0},h_T
\mid \kappa^n,
\nu)-\varepsilon)\\[5pt]
\displaystyle \;\;\;\;\;\;\;\;\;\;\;\;\;\;\;\;\;\;\;\;
\;\;\;\;\;\;\;\;\;\;\;\;\;\;\;\; \;=
 V_{i,0} ^n
\left(x_{i,0},s_0 \mid \kappa^n, {\nu}   \right) -\varepsilon,\;\;\;
\;\;\;\;\;\; \;\;\;\;\;\; \;\;\; \;\;\;\;\;\; \;\;\; \forall n \ge
N, \;\;\;\;\;\; \; \forall x_{i,0} \in \mathcal {X}_0.
\end{array}
\end{equation}
The DOE strategy $\nu$, by definition, maximizes consumer $i$'s
oblivious value function among all possible dynamic oblivious
strategies. The first inequality in \eqref{equa:L_255} follows from
the fact that $\widetilde V_{i,0} ( x_{i,0}, s_0 \mid \kappa^n,
{\nu} )$ is a weighted sum of the oblivious values achieved by a
family of dynamic oblivious strategies\footnote{Note that for a
given $f^n_{-i,0}$, the action taken by $\kappa^n$ depends only on
$x_{i,0}$ and $h_t$, and that $\widetilde V_{i,0} \left(x_{i,0},h_T
\mid \kappa^n, {\nu} \right)$ is the oblivious value achieved by a
mixed strategy; cf. the footnote associated with
\eqref{equa:L_ob}.}. The second inequality in \eqref{equa:L_255} is
due to \eqref{equa:L_21}, and the last equality in
\eqref{equa:L_255} follows from \eqref{equa:L_24}. The desired
result, \eqref{equa:dynamic_mu}, follows.

\textbf{Step 4:} {\it Consumer $i$'s optimal oblivious value can be
asymptotically achieved at an $n$-consumer game under a DOE
strategy.}

In this step, we consider the case where all consumers in an
$n$-consumer model use a DOE strategy $\nu$, and show that
\begin{equation}\label{equa:dynamic_nu}
\mathop {\lim }  \limits_{n \to \infty } \left( \widetilde V_{i,0} (
x_{i,0}, s_0 \mid \nu ,{\nu} ) - V_{i,0} ^n \left( x_{i,0}, s_0,
 \mid \nu, {\nu} \right) \right) =
0, \;\;\; \forall s_0 \in \mathcal {S}, \;\;\; \forall x_{i,0} \in
\mathcal {X}_0.
\end{equation}
According to \eqref{equa:L_21}, with $\kappa^n=\nu$, for any
$\varepsilon>0$, we can find some $N$ such that
$$
\begin{array}{l}
\left| \widetilde V_{i,0} ( x_{i,0},h_T \mid \nu, \nu ) -
V^n_{i,0}(x_{i,0},h_T \mid \nu, \nu)  \right| \le \varepsilon,
\;\;\;\forall n \ge N, \;\;\;\forall h_T \in \mathcal {H}_T,\;\;\;
\forall x_{i,0} \in \mathcal {X}_0.
 \end{array}
$$
The desired result in \eqref{equa:dynamic_nu} follows from
\eqref{equa:L_24} and \eqref{equa:L_224}. Theorem \ref{thm:APE}
follows from \eqref{equa:dynamic_mu} and \eqref{equa:dynamic_nu}.

\section{Proof of Theorem \ref{thm:SO_social}}\label{sec:C}

\subsection{Proof of Part (a)}

 We will show that in a continuum model, a DOE strategy maximizes the
expected social welfare among all possible symmetric dynamic oblivious
strategy profiles (part (a) of the theorem), i.e., that if $\nu$ is DOE, then
\begin{equation}\label{equa:SO_L1}
  \widetilde {\mathcal W}_0(s_0 \mid {\nu})  = \sup\nolimits_{ \vartheta \in \mathfrak{V}}
   \widetilde {\mathcal W}_0(s_0 \mid \vartheta)  ,
\qquad \forall s_0 \in \mathcal {S}.
\end{equation}
Let $S$ and $X$ be the cardinality of $\mathcal {S}$ and $\mathcal
{X}_{0}$, respectively. Given the initial global state $s_0$, the
number of possible histories of length $t+1$ is $S^{t}$. Hence, the
number of all possible histories commencing at state $s_0$ is $
\sum\nolimits_{t=0}^T S^{t}.$ Given an initial global state $s_0$,
the expected social welfare defined in \eqref{equa:SO_CWel} is a
deterministic function of the following $(X\sum\nolimits_{t=0}^T
S^{t})$-dimensional action vector:
\begin{equation}\label{equa:SO_vec}
\left\{\nu_t\left(x, h_t \right) \right\}_{x\in\mathcal {X}_0,\;\;
h_t\in \mathcal {H}(s_0)},
\end{equation}
where $\mathcal {H}(s_0)$ is the set of positive probability
histories commencing at state $s_0$.
 Under Assumption \ref{A:convex}, the expected social welfare defined in
\eqref{equa:SO_CWel} is a concave function of the vector in
\eqref{equa:SO_vec}. Therefore, the following conditions are
necessary and sufficient for the action vector (in the form of
\eqref{equa:SO_vec}) associated with the DOE strategy $\nu$ to
maximize the expected social welfare, among all possible dynamic
oblivious strategies\footnote{We use the notations $\partial_+f$ and
$\partial_-f$ to denote the right and left, respectively,
derivatives of a function $f$.}:
 \begin{equation}\label{equa:SO_condition2}
\begin{array}{l}
\left\{ \begin{array}{l} \displaystyle \dfrac{\partial_{+}
U_t(l_{\nu,h_t}(x_{i,0}), s_t, \nu_t(x_{i,0}, h_t))}{\partial_{+}
\nu_t(x_{i,0}, h_t)} \le \widetilde p_{t|\nu,h_t}+ \widetilde
w_{t|\nu,h_t}
+1_{\tau < T} \cdot g^+_{t|\nu,h_t}(x_{i,0}), \\
\displaystyle\;\;\;\;\;\;\;\;\;\;\;\;\;\;\;\;\;\;\;\;\;\;\;\;\;\;\;\;\;\;\;\;\;\;\;\;\;\;\;\;\;\;\;\;\;\;\;\;\;
\;\;\;\;\;\;\;\;\;\;\;\;\;\;\;\;\;\;\;\;\;\;\;\;\;\;\;\;\;\;\;\;\;\;\;\;\;\;\;\;\;\;\;\;\;\;
\;\;\;{\rm{if}}\;\;\nu_t(x_{i,0}, h_t)
<B,\\
\displaystyle \dfrac{\partial_{-}   U_t(l_{\nu,h_t}(x_{i,0}), s_t,
\nu_t(x_{i,0}, h_t))}{\partial_{-} \nu_t(x_{i,0}, h_t)} \ge
\widetilde p_{t|\nu,h_t}+ \widetilde w_{t|\nu,h_t} +1_{\tau < T}
\cdot g^-_{t|\nu,h_t}(x_{i,0}),
 \end{array} \right.\\
\;\;\;\;\;\;\;\;\;\;\;\;\;\;\;\;\;\;\;\;\;\;\;\;\;\;\;\;\;\;\;\;\;\;\;\;\;\;\;\;\;\;\;\;\;\;\;\;\;
\;\;\;\;\;\;\;\;\;\;\;\;\;\;\;\;\;\;\;\;\;\;\;\;\;\;\;\;\;\;\;\;\;\;\;\;\;\;\;\;\;\;\;\;\;\;\;\;\;\;\;
\;\;\;{\rm{if}}\;\;\nu_t(x_{i,0}, h_t)
>0,
  \end{array}
\end{equation}
where $l_{\nu,h_t}(x_{i,0})$ is consumer $i$'s state, $x_{i,t}$,
under a (positive probability) history $h_t$ and the strategy $\nu$ (cf. p.15), the prices,
$\widetilde p_{t|\nu,h_t}$ and $\widetilde w_{t|\nu,h_t}$
 are given in (\ref{equa:OE_price0}) and
(\ref{equa:OE_price}), and where, if $k_{h_\tau}(\cdot)$ (cf. the
definition in \eqref{equa:SO_k}) is nondecreasing in $a_{i,t}$ for
any $t<\tau \le T$, then $g^+_{t|\nu,h_t}(x_{i,0})$ is given
by\footnote{If for some $\tau>t$, $k_{h_\tau}(\cdot)$ is decreasing
in $a_{i,t}$, then the right partial  derivative of
$U_\tau(x_{i,\tau}, s_\tau, a_{i,\tau})$ with respect to
$z_{i,\tau}$ in \eqref{equa:SO_g} should be replaced by its left
partial derivative.}
 \begin{equation}\label{equa:SO_g}
g^+_{t|\nu,h_t}(x_{i,0})=\E\left\{ \widetilde q_{t+1|\nu,h_{t+1}}
-\sum\limits_{\tau=t+1}^T
 \dfrac{\partial_{+} U_\tau(x_{i,0}, z_{i,\tau}, s_\tau,
a_{i,\tau})}{\partial_{+} z_{i,\tau}} \cdot \dfrac{\partial_{+}
z_{i,\tau}}{\partial_{+} a_{i,t}}\right\},\qquad \forall x_{i,0} \in
\mathcal{X}_0,
\end{equation}
where the price, $\widetilde q_{t+1|\nu,h_{t+1}}$, is defined in
(\ref{equa:OE_price}), the expectation is over the future global
states, $\{s_\tau\}_{t+1}^T$,
$z_{i,\tau}=k_{h_\tau}(x_{i,0},a_{i,0},\ldots,a_{i,\tau-1})$ for
$\tau>t$, and $a_{i,\tau}=\nu_\tau(x_{i,\tau},h_\tau)$ for $\tau \ge
t$. The expression \eqref{equa:SO_g} is the part of the right
derivative of the expected  social welfare \eqref{equa:SO_CWel} with
respect to the action $a_{i,t}$, which reflects the influence of
consumer $i$'s action at stage $t$ on the ancillary cost $\widetilde
H(\widetilde A_t, \widetilde A_{t+1}, \overline s_{t+1})$  at the
next stage,
 and on her future utility (due to the influence of the action $a_{i,t}$ on the future state $z_{i, \tau}$,
  through the functions $k_{h_\tau}(\cdot)$).
In \eqref{equa:SO_condition2}, $g^-_{t|\nu,h_t}(x_{i,0})$ can be
defined by replacing the right (left) partial  derivatives in
\eqref{equa:SO_g} with left (respectively, right) partial
derivatives.

Given an initial global state $s_0$, and the initial state of
consumer $i$, $x_{i,0}$, her oblivious value, defined in
\eqref{equa:OE_Cvalue}, is a deterministic, concave function of the
vector
\begin{equation}\label{eq:vec}
 \left\{\nu_t\left(x_{i,0},
h_t \right) \right\}_{h_t\in\mathcal {H}(s_0)}
\end{equation}
of actions that she would take at any given stage and for any given
history. Since the DOE strategy $\nu$ maximizes consumer $i$'s
oblivious value, it is easily checked that the vector in
\eqref{eq:vec}
 must satisfy the conditions \eqref{equa:SO_condition2}. Since this is true for any $x_{i,0} \in \mathcal
{X}_0$, we conclude that the action vector \eqref{equa:SO_vec}
(which is comprised by putting together the vectors in
\eqref{eq:vec}, for different types in the set $\mathcal{X}_0$)
satisfies the conditions \eqref{equa:SO_condition2}. Thus, the DOE
$\nu$ satisfies the sufficient condition for optimality and the
 result \eqref{equa:SO_L1} follows.

\subsection{Proof of Part (b)}

For a given initial global state $s_0$,
let us fix some initial distribution $f^n_0$ with $d(f^n_0, \eta_{s_0} ) \le \delta$, where
$\delta$ is small. In Step 1, we compare the social welfare achieved by various strategy profiles and show that
$$
  \frac{1}{n}  \mathcal {W}^n_0(f^n_{0}, s_0 \mid \v{\kappa}^n)  \le
   \frac{1}{n}  {\mathcal W}^n_0( s_0 \mid
   \v{\vartheta}^{n,f^n_0}) \approx  \widetilde {\mathcal W}_0( s_0 \mid
   \vartheta^{n,f^n_0}).
$$
Here, $\kappa^n$ is a general history-dependent strategy profile for the $n$-consumer
model (cf. \eqref{eq:his}). The symmetric strategy profile $\v{\vartheta}^{n,f^n_0}=(\vartheta^{n,f^n_0},\ldots,\vartheta^{n,f^n_0})$
is one that maximizes expected
social welfare given the initial population state $f^n_0$. In Step 1, we will
argue that $\vartheta^{n,f^n_0}$ can be identified with a dynamic oblivious strategy. In the approximate equality we are comparing the expected
(over future global states, $\{s_t\}_{t=1}^T$) social welfare under the same oblivious strategy $\vartheta^{n,f^n_0}$
(hence the same sequence of actions for each consumer type $x \in \mathcal{X}_0$) under two different initial population states
(initial distributions
of consumer types), $f^n_0$ and $\eta_{s_0}$.

Since $\nu$ is a DOE, part (a) of the theorem implies that
$$
\widetilde {\mathcal W}_0( s_0 \mid
   \vartheta^{n,f^n_0}) \le \widetilde {\mathcal W}_0( s_0 \mid
   \nu).
$$
Note that as the number of consumer grows large, with high probability
the initial population state $f^n_0$ is close to its expectation,
$\eta_{s_0}$.
In Step 2, we complete the proof of part (b) by showing that
$$
\widetilde {\mathcal W}_0( s_0 \mid
   \nu) \approx \frac{1}{n}  \E \{ {\mathcal W}^n_0(f^n_{0}, s_0 \mid
   \v{\nu}^n)    \},
$$
where the expectation is over the initial population state $f^n_{0}$.

\textbf{Step 1:} {\it If the initial population state is close to
its expectation, the optimal social welfare in an $n$-consumer model
can be approximated by the social welfare achieved by a dynamic
oblivious strategy in the corresponding continuum model. }

In this step, we aim to show that in an $n$-consumer model, for any
given initial global state $s_{0}$ and any $\varepsilon>0$, there
exists some $\delta>0$ such that for any initial distribution
$f^n_{0}$ with $d(f^n_0, \eta_{s_0} ) \le \delta$, we can find a
dynamic oblivious strategy \(\vartheta^{n,f^n_0}\) that satisfies
\begin{equation}\label{equa:SO_L2}
 {\mathcal W}^n_0(f^n_{0}, s_0 \mid \v{\kappa}^n)  \le
  n \widetilde {\mathcal W}_0( s_0 \mid \vartheta^{n,f^n_0}) + \varepsilon n,
\end{equation}
for all symmetric
history-dependent strategy profiles, $\v{\kappa}^n = (\kappa^n,\ldots,\kappa^n)$. Given an initial global state
$s_{0}$ and an initial population state $f^n_{0}$, we observe that
the social welfare, ${\mathcal W}^n_0(f^n_{0}, s_0 \mid
\v{\kappa}^n) $, is a deterministic, concave function of the
following vector of consumers' actions under different histories,
\begin{equation}\label{equa:SO_vec2}
\left\{\kappa^n_{t}\left(m_{i,{\kappa}^n,h_t}(x_{i,0}),h_t,f^n_{-i,t}
\right) \right\}_{h_t\in\mathcal {H}(s_0),\;\;x_{i,0}\in\mathcal
{X}_0, \;\; i=1,\ldots,n},
\end{equation}
where $m_{i,{\kappa}^n,h_t}: \mathcal {X}_0
 \to \mathcal {Y}_t$
maps consumer $i$'s initial state into her augmented state at the
history $h_t$, under the strategy profile $\v{\kappa}^n$, and
$f^n_{-i,t}$ is the distribution of other consumers' augmented
states at the history $h_t$, under the strategy profile
$\v{\kappa}^n$. Note that given the initial population state
$f^n_{0}$, the strategy profile $\v{\kappa}^n$, and a history $h_t$,
the augmented state of consumer $i$ at stage $t$ depends only on her
initial state $x_{i,0}$.

 Since
the social welfare ${\mathcal W}^n_0(f^n_{0}, s_0 \mid \v{\kappa}^n)
$ is concave in the action vector in \eqref{equa:SO_vec2},
 there exists a
symmetric solution, $ \vartheta^{n,f^n_0}=\{
\vartheta^{n,f^n_0}_t\}_{t=0}^T$, such that if at any history $h_t
\in \mathcal
 {H}(s_0)$, all consumers with the same initial state
 take the same action according to
 \begin{equation}\label{eq:v}
 a_{i,t}= \vartheta^{n,f^n_0}_t(x_{i,0},h_t),\;\;\; i=1,\ldots,n,
 \end{equation}
then the expected social welfare, ${\mathcal W}^n_0(f^n_{0}, s_0
\mid \v{\kappa}^n) $, is maximized among all possible symmetric
history-dependent strategy profiles\footnote{The fact that the supremum is attained is a consequence of our continuity assumption and the
fact that the various variables of interest can be restricted to be in a compact set.}. In \eqref{eq:v} we have defined
a dynamic oblivious strategy $\vartheta^{n,f^n_0}$ that maximizes
the expected social welfare in the $n$-consumer model, conditional
on the initial global state being $s_0$, and the initial population
state being $f^n_0$. That is, in an $n$-consumer model, for any
given $s_0$ and $f^n_0$, there exists a dynamic oblivious strategy
$\vartheta^{n,f^n_0}$ such that
\begin{equation}\label{eq:v2}
\sup\nolimits_{ {\kappa}^n} {\mathcal W}^n_0(f^n_{0}, s_0 \mid
\v{\kappa}^n)  =
 {\mathcal W}^n_0(f^n_{0}, s_0 \mid \v{\vartheta}^{n,f^n_0}),
\end{equation}
where $\v{\vartheta}^{n,f^n_0} =
(\vartheta^{n,f^n_0},\ldots,\vartheta^{n,f^n_0})$ is the corresponding symmetric
dynamic oblivious strategy profile. To verify \eqref{equa:SO_L2}, it suffices to show that for any $\varepsilon>0$, there exists some
$\delta>0$ such that for any $f^n_0$ with $d(f^n_0,\eta_{s_0}) \le
\delta$,
\begin{equation}\label{equa:SO_L22}
  \left|  {\mathcal W}^n_0(f^n_{0}, s_0 \mid \v{\vartheta}^{n,f^n_0})  -
  n \widetilde {\mathcal W}_0(s_0 \mid \vartheta^{n,f^n_0}) \right| \le  \varepsilon n,
\end{equation}
i.e., if all consumers use the strategy $\vartheta^{n,f^n_0}$, the
difference between the optimal social welfare achieved in an
$n$-consumer model and the social welfare achieved in the
corresponding continuum model can be made arbitrarily small, if the
initial population state is close enough to its expectation,
$\eta_{s_0}$. We next argue that the result in \eqref{equa:SO_L22}
holds for any dynamic oblivious strategy $\vartheta$.


To prove \eqref{equa:SO_L22}, we first upper bound the difference
between the supplier cost in an $n$-consumer model and that in the
corresponding continuum model. Since all cost functions are
Lipschitz continuous (see Eqs. \eqref{equa:Alip} and
\eqref{equa:Alip1}),
 for any $\varepsilon>0$,  there exists some
$\delta_1>0$ such that if
\begin{equation}\label{equa:SO_L2condtion}
\left|A^n_t - n \widetilde A_{t|\vartheta,h_t}\right| \le X \delta_1
Bn,\qquad t=0,\ldots,T,\qquad \forall h_t \in \mathcal {H}(s_0),
\end{equation}
then
\begin{equation}\label{equa:SO_CC}
\left|C^n(A^n_t,s_t)-C^n(n \widetilde
A_{t|\vartheta,h_t},s_t)\right| \le n
\varepsilon/(3T+3),\qquad t=0,\ldots,T, \qquad \forall h_t \in \mathcal
{H}_t(s_0),
\end{equation}
\begin{equation}\label{equa:SO_HH0}
\left|H_0^n(A^n_0,s_0)-H_0^n(n \widetilde
A_{0|\vartheta,h_0},s_0)\right| \le n \varepsilon/(3T+3),
\end{equation}
and for $t=1,\ldots,T$,
\begin{equation}\label{equa:SO_HH}
\;\;\;\left|H^n(A^n_{t-1},A^n_{t},\overline s_t)-H^n(n \widetilde
A_{t-1|\vartheta,h_{t-1}}, n \widetilde A_{t|\vartheta,h_t},
\overline s_t)\right| \le n \varepsilon/(3T+3), \qquad \forall h_t
\in \mathcal {H}_t(s_0),
\end{equation}
where $\mathcal {H}_t(s_0)$ is the set of all histories of length
 $t+1$ commencing at state $s_0$. Given an initial population state
$f^n_0$, if all consumers use the strategy $\vartheta$, the
aggregate demand under a history $h_t$ is
$$
A^n_t=n\sum\limits_{x \in \mathcal {X}_0} f^n_0(x) \vartheta_t(x,
h_t).
$$
From \eqref{equa:OE_demand} we observe that if $d(f^n_0,\eta_{s_0})
\le \delta_1$,  the condition in \eqref{equa:SO_L2condtion} holds,
and then Eqs. \eqref{equa:SO_CC}-\eqref{equa:SO_HH} are verified.

We now show that if the initial population state is close to its
expectation, the total utility obtained by all consumers is
close to its counterpart in the corresponding continuum model. Given
an initial population state $f^n_0$, we write the total utility
obtained by all consumers under a history $h_t$ as
$$
\sum\limits_{i=1}^n U_t(x_{i,t},s_t,a_{i,t}) =n\sum\limits_{x \in
\mathcal {X}_0} f^n_0(x) U_t\left(
l_{\vartheta,h_t}(x),\vartheta_t(x, h_t) ,s_t\right).
$$
On the other hand, the utility achieved in the corresponding
continuum model is given by
$$
\widetilde U_{t|\vartheta,h_t}  \buildrel \Delta \over =
\sum\limits_{x \in \mathcal {X}_0} \eta_{s_0}(x) U_t\left(
l_{\vartheta,h_t}(x),\vartheta_t(x, h_t) ,s_t\right).
$$
We have that if $d(f^n_0,\eta_{s_0}) \le \varepsilon/(3XQ(T+1))$,
then
\begin{equation}\label{equa:SO_UU}
\left|\sum\nolimits_{i=1}^n U_t(x_{i,t},s_t,a_{i,t})- n\widetilde
U_{t|\vartheta,h_t}  \right| \le n
\varepsilon/(3T+3),\;\;\;t=0,\ldots,T,\;\;\;\forall h_t \in \mathcal
{H}_t(s_0),\;\;\;\forall n \in \mathbb{N}^+,
\end{equation}
Let $\delta=\min\{\delta_1, \varepsilon/(3XQ(T+1))\}$. If
$d(f^n_0,\eta_{s_0}) \le \delta$, from
\eqref{equa:SO_CC}-\eqref{equa:SO_UU} we have
$$
\left| \widetilde W_t(h_t \mid \vartheta ) -
 {W}^n_{t}(f^n_t,h_t \mid \v{\vartheta}^n) \right| \le n
\varepsilon/(T+1),\;\;\;t=0,\ldots,T,\;\;\;\forall h_t \in \mathcal
{H}_t(s_0).
$$
Eq. \eqref{equa:SO_L22} follows from the definition of expected
social welfare in an $n$-consumer model
 \eqref{equa:SO_welfarenu}, and in a continuum model \eqref{equa:SO_CWel}.
 The desired result in \eqref{equa:SO_L2} follows.

\bigskip


\textbf{Step 2:} {\it Asymptotic social optimality of a DOE. }

In this step, we complete the proof of part (b) of the theorem, using the fact that as the
number of consumers grows large, with high probability the initial
population state is close to its expectation. Note that the action space is $[0,B]$, so that
$|A_t| \le nB $. Using Assumption \ref{A:continuous}, $C^n(\cdot)/n$ is therefore bounded. A similar argument
holds for $H_0^n(\cdot)/n$ and $H^n(\cdot)/n$. Furthermore, the total utility per consumer is also bounded.
Thus, there exists some constant $D$ that upper bounds
$|\mathcal {W}_0^n/n|$. We define $\mathfrak{F}_{s_0}^n(\delta)$ as the set of initial population
states such that $d(f^n_0, \eta_{s_0}) \le \delta$.
By the law of large numbers, for any pair of positive real numbers, $\varepsilon$ and $\delta$,
we can find an integer $N$ such that
\begin{equation}\label{equa:SO_L32}
\sum\limits_{f_0^n \notin
\mathfrak{F}_{s_0}^n(\delta)} \mathbb{P}\Big(F^{n}_{s_0} = f_0^n
\Big) \cdot \sup_{\kappa^n} \left| \mathcal {W}^n_{0}(f^n_0,s_0 \mid \v{\kappa}^n) \right| \le  D \mathbb{P}( d(f_0^n, \eta_{s_0}) > \delta )  \le \varepsilon
n,\qquad \forall n \ge N.
\end{equation}

For any $\varepsilon>0$, let $\delta$ be the positive real number
defined in \eqref{equa:SO_L2}, and let $N$ be the positive integer
given in \eqref{equa:SO_L32}; for any $n \ge N$ and any symmetric
history-dependent strategy profile $\v{\kappa}^n$,  we have
$$
\begin{array}{l}
\displaystyle \;\;\;\; \E\left\{  \mathcal {W}^n_{0}(f^n_0,s_0 \mid
\v{\kappa}^n)  \right\}\\[4pt]
\displaystyle \le \sum\nolimits_{f_0^n \in
\mathfrak{F}_{s_0}^n(\delta)} \mathbb{P}\big(F^{n}_{s_0} = f_0^n
\big) \cdot \mathcal {W}^n_{0}(f^n_0,s_0 \mid \v{\kappa}^n)
+ \varepsilon n\\
 \displaystyle \le   \sum\nolimits_{f_0^n \in
\mathfrak{F}_{s_0}^n(\delta)} \mathbb{P}\big(F^{n}_{s_0} = f_0^n
\big) \cdot \left( n \widetilde {\mathcal W}_0(s_0 \mid
\vartheta^{n,f_0^n}) + \varepsilon n   \right)+ \varepsilon n\\
 \displaystyle \le   \sum\nolimits_{f_0^n \in
\mathfrak{F}_{s_0}^n(\delta)} \mathbb{P}\big(F^{n}_{s_0} = f_0^n
\big) \cdot \left(  n \widetilde {\mathcal W}_0(s_0 \mid \nu) +
\varepsilon n \right)+ \varepsilon n\\
 \displaystyle \le   \sum\nolimits_{f_0^n \in
\mathfrak{F}_{s_0}^n(\delta)} \mathbb{P}\big(F^{n}_{s_0} = f_0^n
\big) \cdot \left( {\mathcal W}^n_0(f^n_{0}, s_0 \mid \v{\nu}^n)  +
2\varepsilon n
\right)+
\varepsilon n\\[7pt]
\displaystyle \le \E \left\{ \mathcal {W}^n_{0}(f^n_0,s_0 \mid
\v{\nu}^n)  \right\}  + 4 \varepsilon n,
 \end{array}
$$
where the first inequality follows from
\eqref{equa:SO_L32}, the second inequality is due to
\eqref{equa:SO_L2}, the third inequality follows from the optimality property of the DOE $\nu$
(part (a) of the theorem),
 the fourth inequality follows similar to
\eqref{equa:SO_L22} (the proof of Eq. \eqref{equa:SO_L22}
remains valid for any dynamic oblivious strategy), and the last inequality follows from
\eqref{equa:SO_L32}.

\section{Numerical Results}\label{sec:num}

In this section we give a numerical example to compare the proposed
pricing mechanism with marginal cost pricing. {The comparison is
carried out in terms of DOEs and the resulting social welfare under
the corresponding continuum model. Towards this purpose, we first
define the DOE for a continuum model under the marginal cost pricing
mechanism,} in Section \ref{sec:num1}. In Section \ref{sec:num2}, we
consider a two-stage dynamic model in which {the} consumers'
marginal utility and
 demand increase at the second stage. We calculate the
equilibria resulting from the two pricing mechanisms, and compare
 the potential of the two pricing mechanisms to improve
social welfare and reduce peak load.

\subsection{Equilibrium under Marginal Cost Pricing}\label{sec:num1}
In an $n$-consumer model, at stage $t \ge 1$, the supplier's
marginal cost is
\begin{equation}\label{eq:equ_mar}
(C^n)'(A^n_t,s_t)+ \dfrac{\partial H^n(A^n_{t-1}, A^n_t, \overline
s_t )}{\partial A^n_{t}} = p^n_t + w^n_t,\;\;\;t=1,\ldots,T.
\end{equation}
At stage $0$, the supplier's marginal cost is
\begin{equation}\label{eq:equ_mar0}
(C^n)'(A^n_0,s_0)+ (H^n_0)'(A^n_0,s_0) = p^n_0 + w^n_0.
\end{equation}
Under marginal cost pricing, each consumer's stage payoff is
\begin{equation}\label{equa:E_Spayoff}
\pi(y_{i,t},\overline s_t,a_{i,t},f^n_{-i,t},u^n_{-i,t}) =
U(x_{i,t},s_t,a_{i,t}) - (p^n_t + w^n_t) \cdot a_{i,t},
\end{equation}
where the stage marginal cost, $p^n_t+w^n_t$, is given in
\eqref{eq:equ_mar} and \eqref{eq:equ_mar0}, and
$y_{i,t}=(x_{i,t},a_{i,t-1})$.

For marginal cost pricing, we now define the nonatomic equilibrium
concept in the corresponding continuum model.  Suppose that all
consumers other than $i$ use a dynamic oblivious strategy $\nu$.
Consumer $i$'s oblivious stage value under  marginal cost pricing
 is given by
 \begin{equation}\label{equa:ex_Cpayoffa}
\widetilde \pi_{i,t}(y_{i,t},\overline s_t , f_{t|\nu, h_t},a_{i,t}
\mid
 {\nu} )   =
 U_t(x_{i,t},s_{t}, a_{i,t}) - (\widetilde p_{t|\nu,h_t} +  \widetilde w_{t|\nu,h_t} ) \cdot a_{i,t},
\end{equation}
where $\widetilde p_{t|\nu,h_t} $ and $\widetilde w_{t|\nu,h_t}$ are
defined in (\ref{equa:OE_price0}) and (\ref{equa:OE_price}).
Replacing the oblivious stage value function in
(\ref{equa:OE_Cpayoffa}) with that given in
(\ref{equa:ex_Cpayoffa}), we can define an equilibrium concept for
 the marginal cost pricing mechanism in a similar way as for the DOE in Section \ref{sec:DOE}.

\subsection{Numerical Example}\label{sec:num2}

In current  wholesale electricity markets, we observe that the
highest daily wholesale price usually occurs when the system load
increases quickly (cf. Fig.\ \ref{Fig:211} in Section
\ref{sec:intro}). Inspired by the above observation, we construct a
two-stage dynamic model, in which the aggregate demand increases
quickly at the second stage,
 to compare the performance of the proposed mechanism with marginal cost pricing.
 For simplicity, we assume that there is a continuum of identical
consumers indexed by $i \in [0,1]$. Each consumer would like to
consume $1+x$ and $1.2-x$ at the two stages, where $x \in [0,E]$.
Here, $E \in [0,0.1]$ (a given constant) is the amount of
electricity demand that can be shifted from the second stage to the
first stage. The value of $E$ will be called demand
substitutability\footnote{
 There are two types
of elasticity of consumers' demand: (i) consumers may curtail their
demand at a high price, and (ii) they may shift their demand to a
less expensive time. The first type of demand response is a price
elasticity, and the second type is an elasticity of substitution
across time. The first type of elasticity is incorporated in our
model through the utility functions, and the second type of
elasticity is incorporated through $E$.}.

Formally, consumer $i$'s state at each stage denotes the maximum
amount of electricity she could use at the stage\footnote{Since all
consumers are {of} the same type, the consumer state space in this
example is a subset of {$[0,\infty)$.}}. For a given consumer $i$,
we have
 $x_{i,0}=1+E$, and her state at stage $1$ is determined as follows:
\begin{enumerate}
\item if $a_{i,0} \le 1$, the maximum
amount of electricity she could use at stage $1$  is $1.2-E+E$,
i.e., $x_{i,1}=1.2$;

\item if $ 1 < a_{i,0} \le x_{i,0}$, the maximum
amount of electricity she could use at stage $1$ is ${x_{i,1}=}
1.2-E+(x_{i,0}- a_{i,0})=2.2-a_{i,0}$;

\item if $ x_{i,0} < a_{i,0}$, the maximum
amount of electricity she could use at stage $1$ is
${x_{i,1}=}1.2-E$.
\end{enumerate}
To summarize, we have
$$
x_{i,1}=1.2-E+\max\{0, x_{i,0}- \max\{ a_{i,0}, 1 \} \}.
$$
 For
each stage $t$, the utility functions are given by
$$
U_t(x_{i,t},s_t,a_{i,t}) = \left\{ \begin{array}{l}
 d_t a_{i,t},\;\;\;\;\;\;\;\;\;\; \;\;\;\;\;\;\; {\rm{if}}\;\;\;0 \le a_{i,t} \le x_{i,t}, \\
 d_t x_{i,t} ,\;\;\;\;\;\;\;\;\;\;\;\;\;\;\;\;\;{\rm{if}}\;\;\;\; \; a_{i,t} > x_{i,t}, \\
 \end{array} \right.
$$
where the slopes are $d_0=10$ and $d_1=12$. Here, we assumed that
the consumers place a larger value on electricity during peak hours,
and that shifting peak load to off-peak hours hurts consumer
utility. For example, rescheduling kitchen and laundry activities
may cause inconvenience for residential consumers; similarly,
industrial consumers may face higher labor cost premiums for
off-peak production.

The primary cost function (cf. Section \ref{sec:model}) is
$\widetilde C(A,s)=A^2$, for any $s$. We assume that the capacity
available at each stage is proportional to the system load, i.e.,
$$
G_t=b_t A_t,\;\;\; t=0,1,
$$
and that the ancillary cost depends only on the difference between
the capacity available at two consecutive stages. At the second
stage (peak hour), we assume that the system operator maintains a
reserve margin of $10\%$, i.e., $b_1=1.1$. We will consider two
different system operator policies: (i)  the system operator does
not forecast the load jump at the second stage, and uses a
conservative policy under which $b_0=1.12$, and (ii)  the system
operator predicts the load jump at the second stage, and ramps up
the system capacity in advance, by letting $b_0=1.2$.

For simplicity, we use a quadratic function to approximate the
ancillary cost associated with load fluctuations:
$$
\widetilde H_0 (A_0,s_0)= 10 (\max\{b_0 A_0 - 1.12, 0
\})^2,\qquad \widetilde H (A_0,A_1,\overline s_1) =20 (\max\{b_1 A_1
- b_0A_0,0 \})^2,
$$
 where $1.12$ represents the capacity available at the stage before the
initial stage\footnote{Suppose that the load at stage ``$-1$'' is
$1$, and that the capacity available at stage $-1$ is $1.12$, under
an average reserve margin of $12\%.$}. We assumed a higher
coefficient, $20$,
 for the ancillary cost at the second stage, due to the
increase of the system load.

For different levels of demand substitutability $E$,
 and two different system operator policies ($b_0$ equal to $1.12$ or $b_0=1.2$), we compare the social welfare (in Section \ref{sec:num21})
and the peak load (in Section \ref{sec:num22}) resulting from the
equilibria of the two pricing mechanisms.

\subsubsection{Social welfare gain.}\label{sec:num21}

For various levels of demand substitutability ($E \in [0,0.1]$), and
the two different system operator policies, we calculate the
equilibria resulting from the two pricing mechanisms. Fig.\ 3
compares the social welfare achieved by the proposed mechanism and
the marginal cost pricing mechanism. We observe from Fig.\ 3 the
following.

\begin{enumerate}{\itemsep=0pt}
\item \textbf{System operator's policy:}
When the consumers have a low level of demand substitutability, the
policy with $b_0=1.2$ achieves a much higher social welfare than the
conservative policy ($b_0=1.12$), under both the proposed and the
marginal cost pricing mechanisms. (This is to be expected, because
when $b_0=1.12$, and with the demand at stage 1 more or less fixed,
the difference $b_1A_1-b_0 A_0$ is necessarily large.) For consumers
with a high level of demand substitutability, the policy with
$b_0=1.2$ achieves a slightly smaller social welfare than the
conservative policy ($b_0=1.12$), because the policy with $b_0=1.2$
results in a lower price at the second stage than the conservative
{one}, and therefore {does not provide enough encouragement to the
consumers to shift} their peak load (cf. the discussion in Section
\ref{sec:num22}).

\item \textbf{Social welfare gain at a low level of demand substitutability:}
At a low level of demand substitutability, e.g., when $E \le 0.02$,
{and under the system operator's conservative policy ($b_0=1.12$),}
we observe that the proposed pricing mechanism achieves
significantly more social welfare gain (the social welfare achieved
by flat rate pricing\footnote{Under flat rate pricing, consumers pay
a fixed (time-invariant) retail price for the electricity they
consume. Since the average retail price is less than {the}
consumers' marginal utility (see Tables \ref{t:price} and
\ref{t:price1}), the payoff-maximizing consumer demand at the two
stages is $1$ and $1.2$, respectively. Since all consumers are
identical, the aggregate demand at the two stages is $1$ and $1.2$.}
is used a reference) than marginal cost pricing; if the system
operator ramps up the capacity in advance ($b_0=1.2$),
 both pricing mechanisms achieve approximately the same social welfare as flat rate pricing.

\item \textbf{Social welfare gain at a high level of demand substitutability:} If the consumers have a high demand substitutability,
e.g., when $E \ge 0.08$, the proposed pricing mechanism achieves
approximately $5\%$ more social welfare gain than marginal cost
pricing under the system operator's conservative policy
($b_0=1.12$); if the system operator ramps up the capacity in
advance ($b_0=1.2$), the proposed pricing mechanism achieves
approximately $50\%$ more social welfare gain than marginal cost
pricing.

\end{enumerate}

 \begin{figure}\label{Fig:welfare}
    \includegraphics[width=10cm]{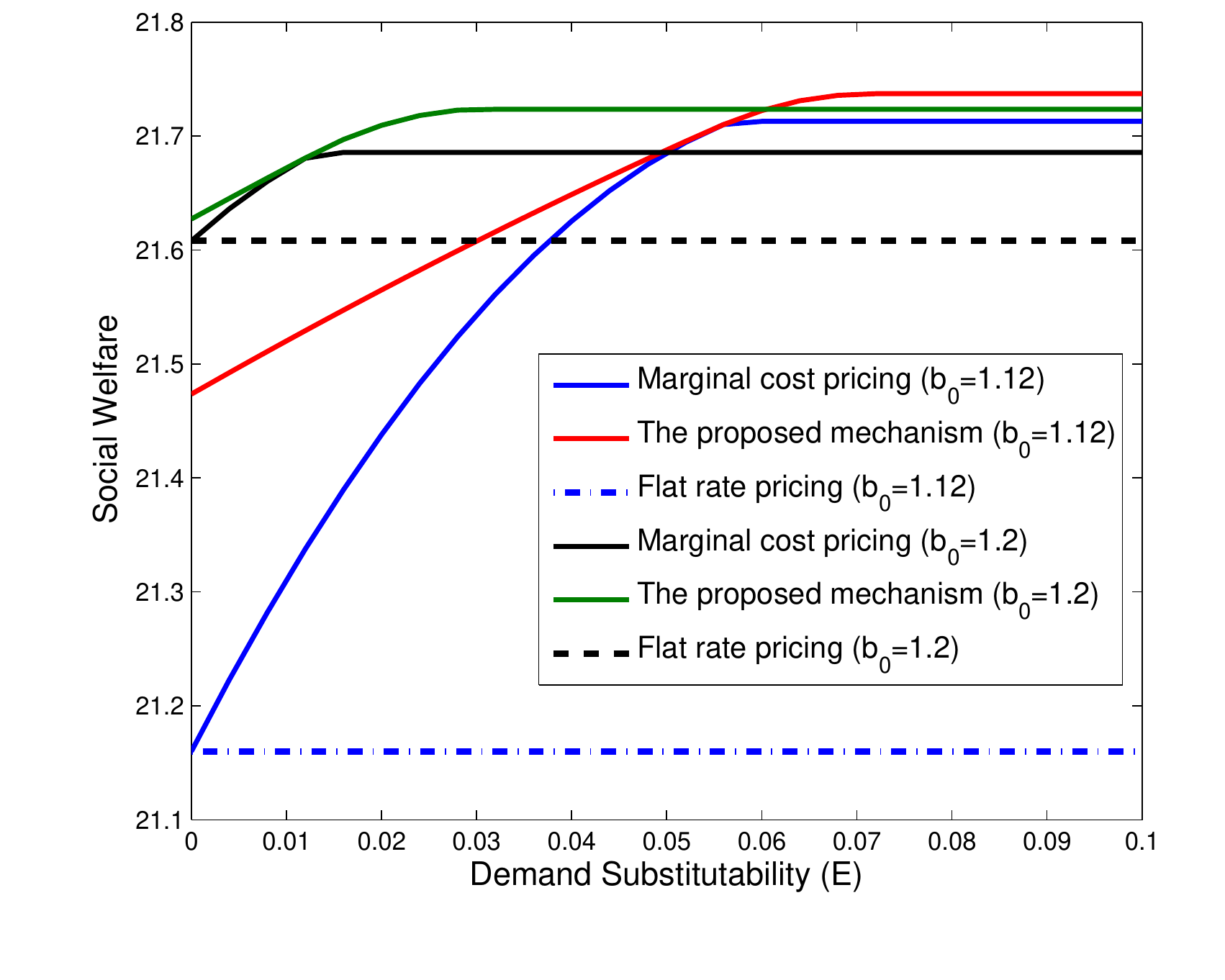}
    \centering
    \caption{The social welfare achieved by the proposed pricing mechanism, the marginal cost pricing mechanism and the flat rate pricing mechanism,
    as a function of the demand substitutability, $E$.}
  \end{figure}

\bigskip
Let us now derive some insights by considering the special case of
zero demand substitutability ($E=0$) and $b_0=1.12$. The one-stage
aggregate demand and the social welfare resulting from the three
pricing mechanisms are given in Table \ref{t:demand}. The prices
faced by consumers are given in Table \ref{t:price}, where the
average retail price is the ratio of the total money a consumer pays
at an equilibrium to her total demand during the two
stages.\footnote{Note that only consumers under flat rate pricing
pay this price. We list the average prices for the two real-time
pricing mechanisms to compare the consumers' expense under different
pricing mechanisms.}
 Note that under the
proposed pricing mechanism, a consumer pays
$$
(p_0+w_0+q_1)a_{i,0}+(p_1+w_1)a_{i,1},
$$
while she would pay $ (p_0+w_0)a_{i,0}+(p_1+w_1)a_{i,1}$ under
marginal cost pricing. {In general, the price $q_1$ will be negative
and will be even smaller if we were to increase} the aggregate
demand at the second stage.  That is, a higher peak load results in
a lower price at the first stage, which encourages consumers to
increase their demand at the
      off-peak hour, {even if they do not derive any additional utility from such an increase. In fact,}
      from Table \ref{t:price} we observe that at the DOE,
the proposed pricing mechanism
      offers each consumer a zero total price on $a_0$.
This may appear illogical at first sight. The reason is that due to
the conservative reserve policy, with $b_0=1.12$, a demand of
$a_0=1$ results in a large increase from $b_0a_0$ to $b_1 a_1$ and
hence a large ancillary cost. The increase of the demand $a_0$
beyond $1$ does not provide any utility to the consumer, but reduces
the ancillary cost. Thus, the counterintuitive choice of
$a_0=1.0901$ serves to mitigate a conservative and somewhat
deficient   reserve policy. {This suggests that further research is
needed that will include an intertemporal optimization of the
reserve policy as well.}

\begin{table} \caption{Demand and social welfare (per
consumer) at $E=0$ and $b_0=1.12$}\label{t:demand}
 \centering
\begin{tabular}{c  c c  c}
 \hline\hline
  & $a_0$ & $a_1$  & Social welfare \\ [0.5ex]
 \hline
 Flat rate & 1 & 1.2  & 21.16\\
Marginal cost & 1 & 1.2  & 21.16\\
Proposed  & 1.0901 &  1.2  & 21.4735\\
 \end{tabular}
\label{table:nonlin}
 \end{table}

  \begin{table} \caption{Price fluctuation at $E=0$ and $b_0=1.12$. The price $p_t+w_t$ equals the
    marginal cost at stage $t$.}\label{t:price}
\centering
\begin{tabular}{c c c c  c}
 \hline\hline
  & $p_0+w_0$ & $p_1+w_1$  & $q_1$ & Average (retail) price\\ [0.5ex]
 \hline
Flat rate  & 2 & 11.2   &- &  7.0182\\
Marginal cost  & 2 & 11.2  &- &  7.0182\\
Proposed & 4.4401 & 6.7609  & -4.4401 &  5.6562\\
 \end{tabular}
\label{table:nonlin}
 \end{table}

For the case where the system operator ramps up the capacity in
advance $(b_0=1.2)$, and consumers have a high level of demand
substitutability ($E=0.08$), the one-stage aggregate demand and the
social welfare resulting from the three pricing mechanisms are given
in Table \ref{t:demand1}. The prices faced by consumers are given in
Table \ref{t:price1}. From Table \ref{t:demand1} we observe that
under the proposed pricing mechanism, consumers would like to shift
$ 0.031$ peak load to off-peak hours, while under the marginal cost
pricing mechanism, consumers are willing to shift less than $0.014$
peak load to off-peak hours. Compared to marginal cost pricing,  the
more flattened load curve resulting from the proposed pricing
mechanism leads to
 $50 \%$ more social welfare gain.

\begin{table} \caption{Demand and social welfare (per
consumer) at $E=0.08$  and $b_0=1.2$}\label{t:demand1}
 \centering
\begin{tabular}{c  c c  c}
 \hline\hline
  & $a_0$ & $a_1$  & Social welfare \\ [0.5ex]
 \hline
 Flat rate & 1 & 1.2  &21.608\\
Marginal cost &  1.0131 &  1.1869  &21.6857\\
Proposed &  1.0308 &   1.1692 &    21.7237\\
 \end{tabular}
 \end{table}

  \begin{table} \caption{Price fluctuation  at $E=0.08$  and $b_0=1.2$. The price $p_t+w_t$ equals the
    marginal cost at stage $t$.}\label{t:price1}
\centering
\begin{tabular}{c c c c c }
 \hline\hline
  & $p_0+w_0$ & $p_1+w_1$ & $q_1$  & Average (retail) price\\ [0.5ex]
 \hline
Flat rate  & 3.92 & 7.68  &- & 5.971\\
Marginal cost  &   4.324 &    6.324  &- &  5.403\\
Proposed  &  4.868 &   4.505 &   -2.363 & 3.568\\
 \end{tabular}
 \end{table}

For a given load curve, the proposed pricing mechanism results in a
larger price difference between stage $1$ and stage $0$ than
marginal cost pricing, because of the negative price $q_1$. The
negative price $q_1$ creates an additional incentive for consumers
to shift their load from stage $1$ to stage $0$. In this way, the
proposed pricing mechanism results in a more flattened load curve
and a higher social welfare than marginal cost pricing (cf. Table
\ref{t:demand1}).

\subsubsection{Peak load reduction.}\label{sec:num22}
Under flat rate pricing, the peak load (the aggregate demand at the
second stage) is $1.2$, because consumers do not have an incentive
to shift their load to off-peak hours.
  Given a pricing
mechanism and a system operator's policy ($b_0$), consumers are
willing to substitute across time only up to a certain level. Even
with a
 high level of demand substitutability,
consumers prefer not to shift much of their peak load, to avoid the
utility loss caused by peak load shifting. For example, with
$b_0=1.2$ and $E=0.08$, consumers  under marginal cost pricing
choose to shift at most $0.013$ peak load (cf. Table
\ref{t:demand1}). In Fig.\ 4, for different values of $b_0$, we
compare the maximum amount of peak load consumers choose to shift
under the proposed pricing mechanism and the marginal cost pricing
mechanism.

\old{

\begin{figure}
\begin{minipage}[t]{0.5\linewidth}
\includegraphics[width=3.03in]{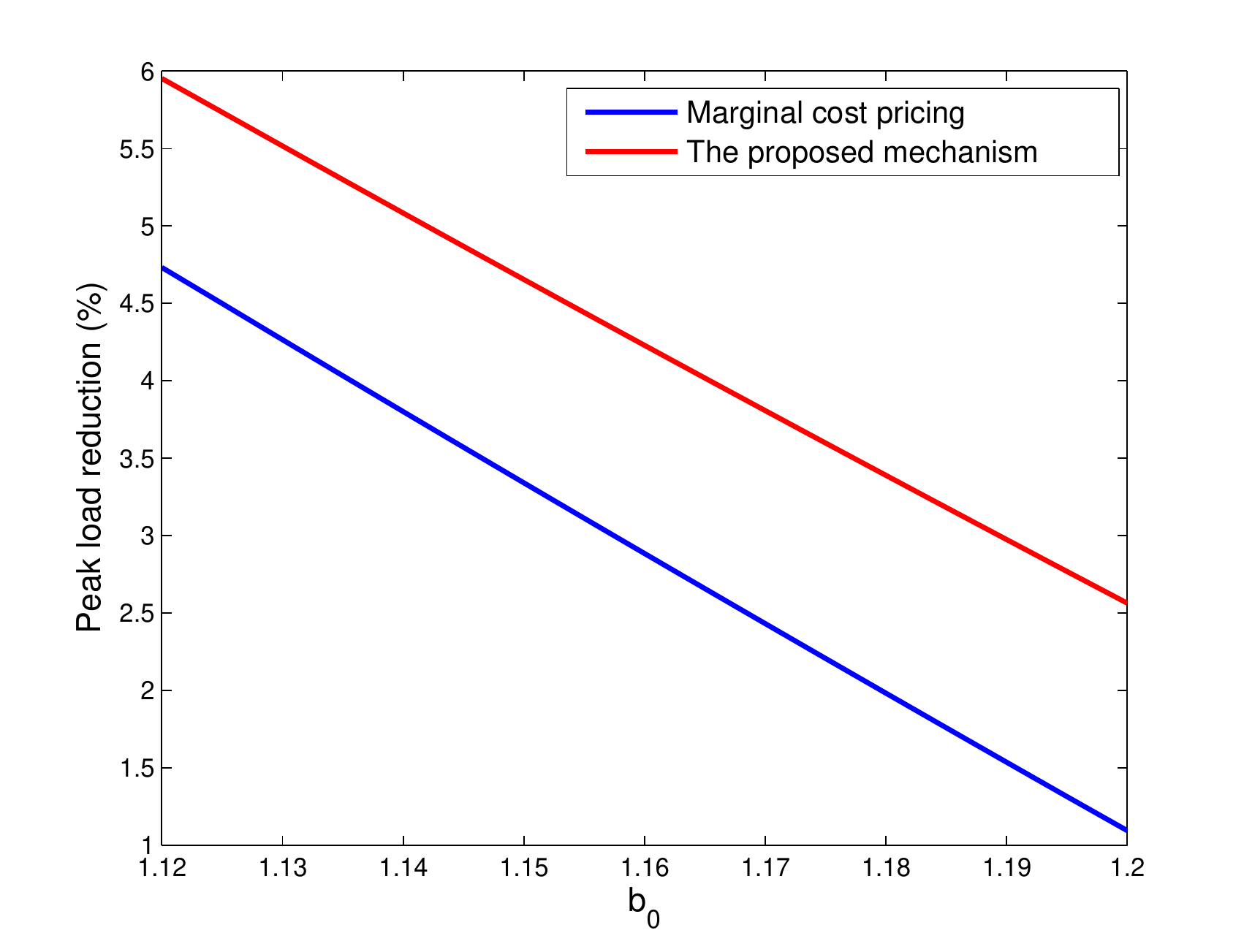}
   \caption{Comparison of the percentage of peak load reduction
    (the peak load under flat rate pricing, $1.2$, is used a reference) resulting from the proposed pricing mechanism and the
    marginal cost pricing mechanism,
    as a function $b_0$.}
\end{minipage}%
\begin{minipage}[t]{0.5\linewidth}
\centering
\includegraphics[width=3.01in]{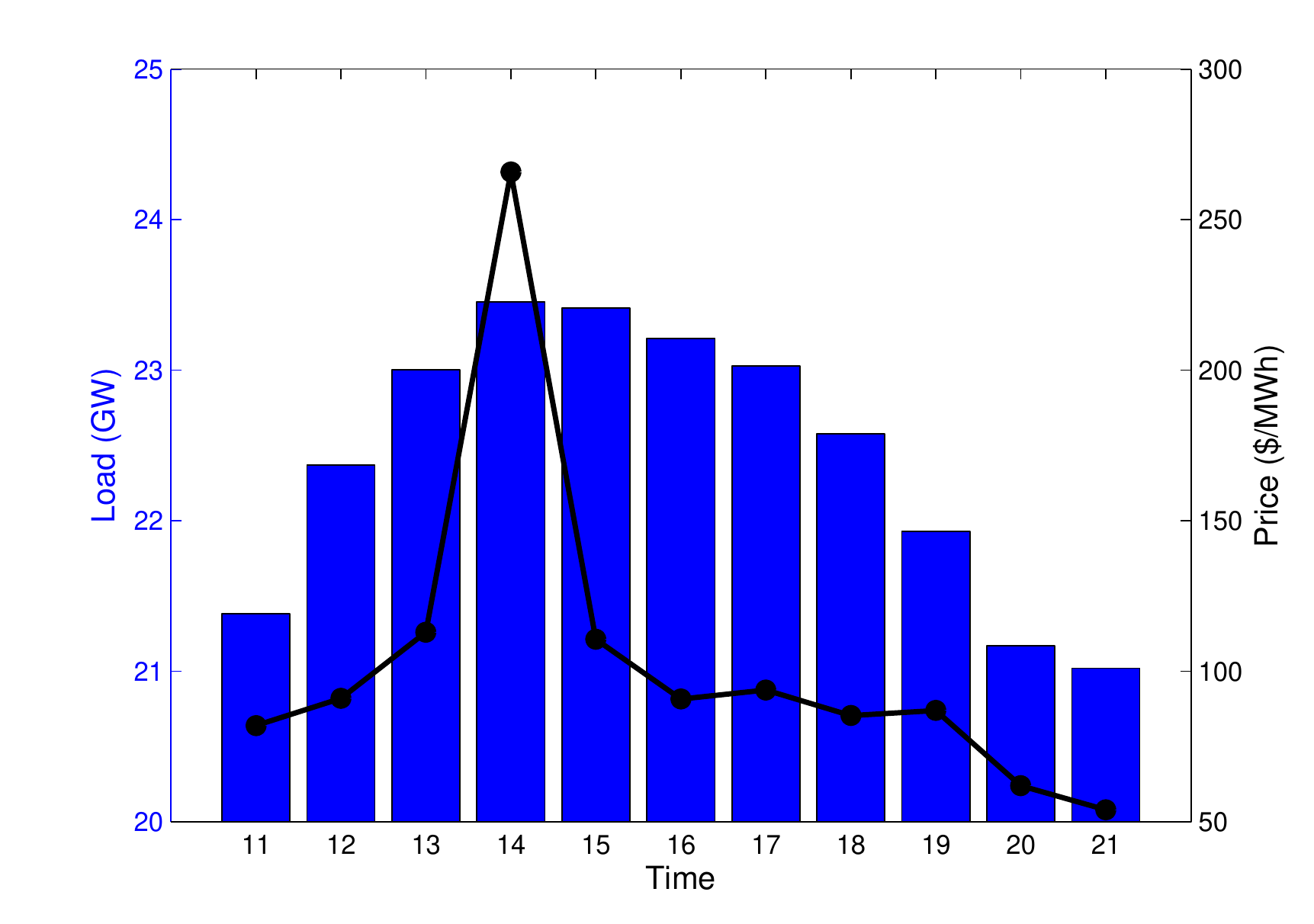}
\end{minipage}
    \caption{Real-time prices and actual system loads on August 01, 2011, ISO
New England Inc. Blue bars represent the real-time system loads and
the dots connected by a black line represent the hourly prices.}
\end{figure}

}

 \begin{figure}
    \includegraphics[width=10cm]{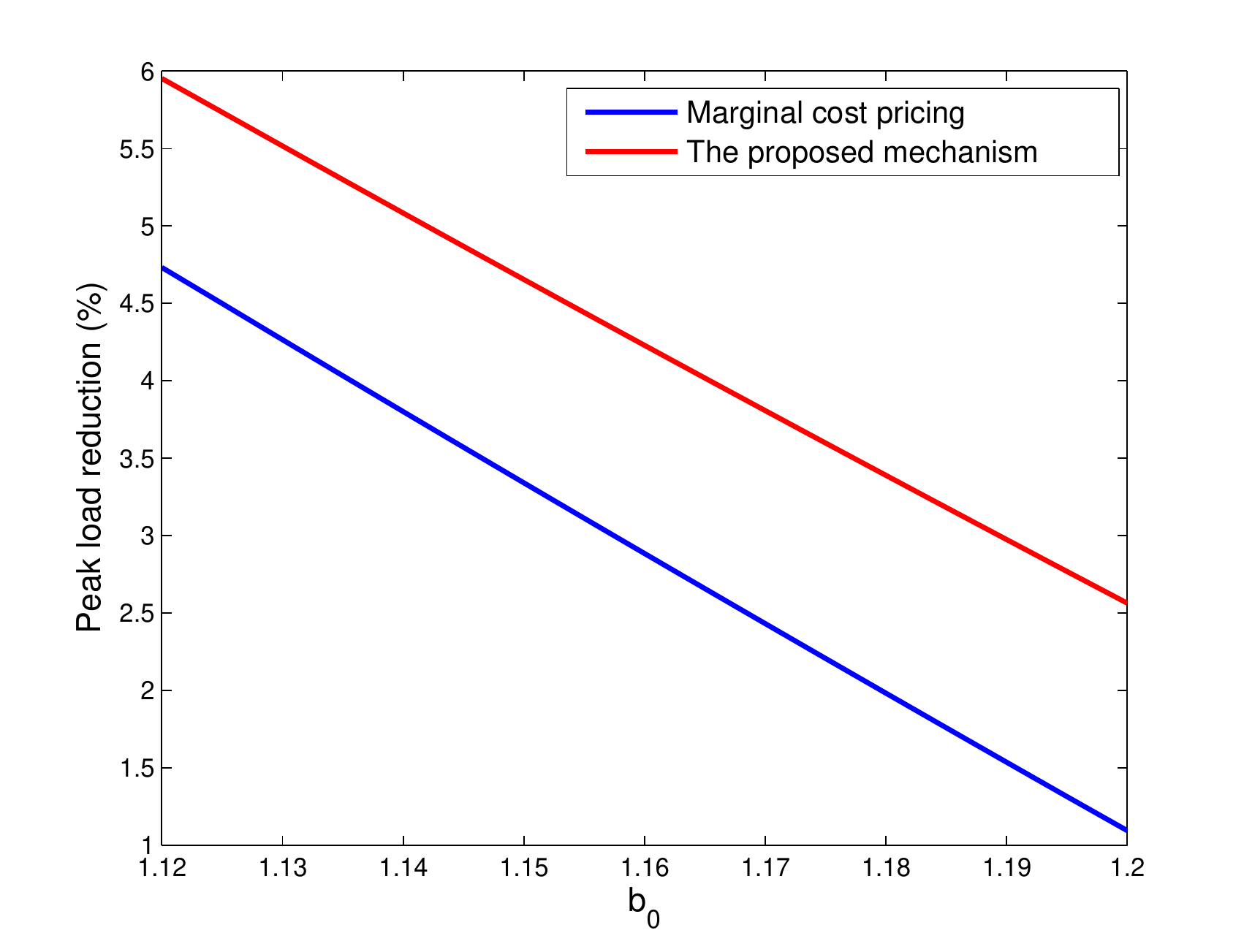}
    \centering
    \caption{Comparison of the percentage of peak load reduction
    (the peak load under flat rate pricing, $1.2$, is used a reference) resulting from the proposed pricing mechanism and the
    marginal cost pricing mechanism,
    as a function $b_0$.}
  \end{figure}

 \begin{figure}
    \includegraphics[width=10cm]{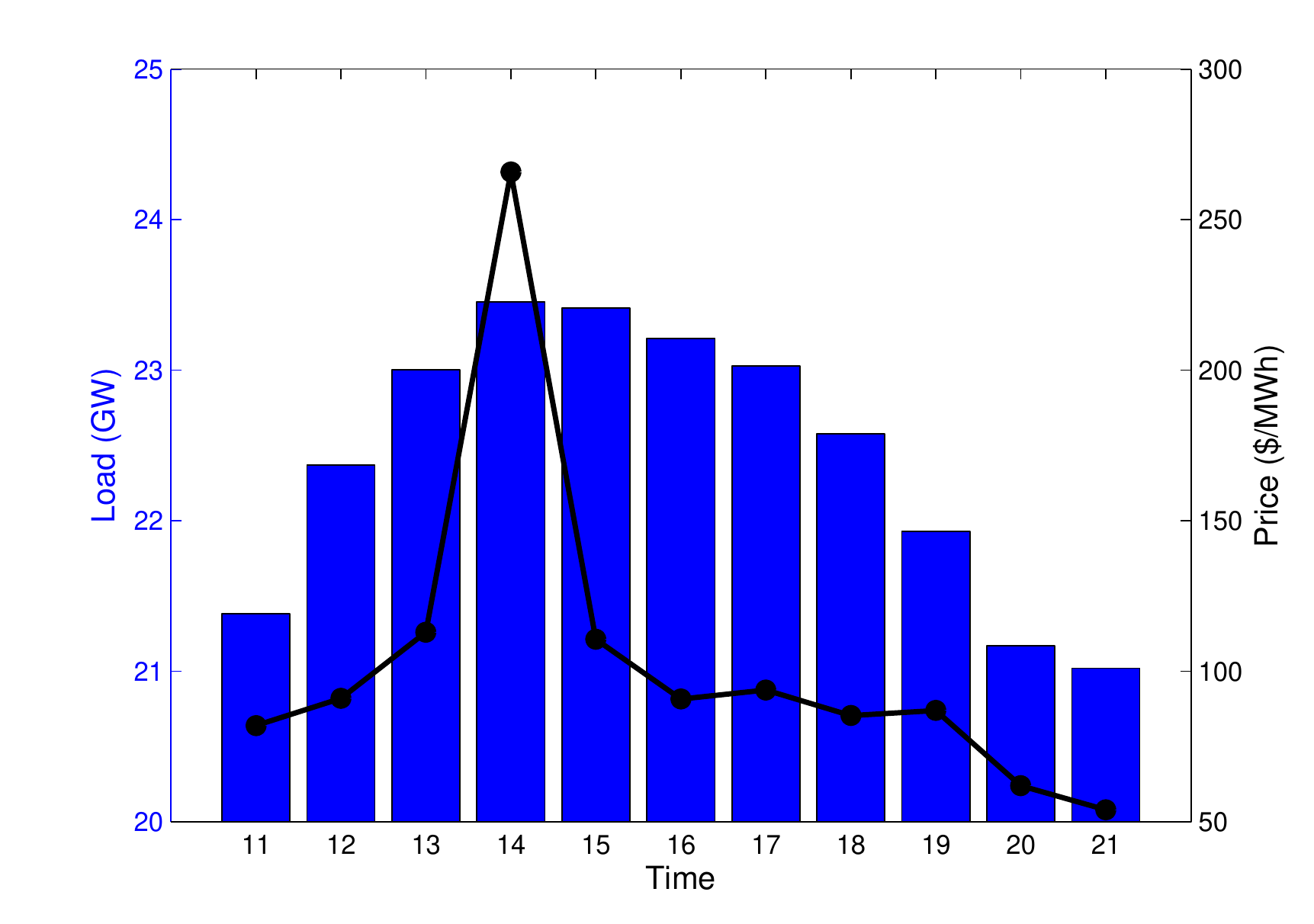}
    \centering
    \caption{Real-time prices and actual system loads on August 01, 2011, ISO
New England Inc. Blue bars represent the real-time system loads and
the dots connected by a black line represent the hourly prices.}
  \end{figure}

 We observe from Fig.\ 4 that the amount of peak load consumers will shift decreases with
 $b_0$. This is
because a larger reserve at the first stage lowers the price at the
second stage, which in turn discourages consumers from shifting
their peak load. The proposed pricing mechanism results in a peak
load which is approximately $1.5$ percent lower than that resulting
from marginal cost pricing, regardless of the value of $b_0$. If the
system operator ramps up the system capacity in advance ($b_0=1.2$),
marginal cost pricing reduces  the system peak load resulting from
flat rate pricing by approximately one percent. Compared to marginal
cost pricing, the negative price $q_1$ in the proposed mechanism
encourages consumers to make a larger shift of their peak load (cf.
the discussion at the end of Section \ref{sec:num21}).

 Fig.\ 5 plots the real-time system loads
and prices on August 1, 2011, a typical hot summer day in New
England\footnote{\url{www.ferc.gov/market-oversight/mkt-electric/new-england/2011/08-2011-elec-isone-dly.pdf}}.
If consumers are able to shift some of their load to the morning
(possibly at the expense of losing some utility), the proposed
pricing mechanism encourages consumers to shift more of their peak
load than marginal cost pricing. Since the highest peak load
determines the generation capacity necessary for system reliability,
the proposed pricing mechanism has a greater potential to reduce the
long-term capacity investment.

\end{appendices}

\end{document}